\documentclass[journal,final,a4paper]{IEEEtran}
\IEEEoverridecommandlockouts

\usepackage{amsmath, amssymb, amsfonts, amsthm, booktabs, cancel, comment, enumitem, eurosym, float, graphicx, mathtools, multirow, nicefrac, placeins, soul, steinmetz, stfloats, siunitx, svg, url, xcolor, tabularx}%
\usepackage[caption=false, font=footnotesize]{subfig}
\usepackage{hyperref}
\usepackage{nomencl}
\usepackage{cleveref}
\makenomenclature

\usepackage{etoolbox}
\renewcommand\nomgroup[1]{%
  \item[
  \ifstrequal{#1}{I}{\textit{Sets and Indices}}{%
  \ifstrequal{#1}{P}{\textit{Parameters}}{%
  \ifstrequal{#1}{V}{\textit{Variables}}{}}}
]}

\usepackage[europeanresistors]{circuitikz}

\floatstyle{ruled}
\newfloat{model}{thp}{lop}
\floatname{model}{\textsc{Model}}

\DeclareMathOperator*{\minimize}{min.}
\DeclareMathOperator*{\maximize}{max.}

\DeclareMathOperator{\Expected}{{\mathbb{E}}}

\newcommand{\AZ}{{\soc{A_i\Sigma^{\nicefrac{1}{2}}}}}
\newcommand{\AZsq}{{\soc{A_i\Sigma^{\nicefrac{1}{2}}}^2}}

\newcommand{\T}{{\!\top}}
\newcommand{\X}{\boldsymbol{X}}
\newcommand{\Y}{\boldsymbol{Y}}
\newcommand{\p}{\boldsymbol{p}}
\newcommand{\om}{\boldsymbol{\omega}}
\newcommand{\Om}{\boldsymbol{\Omega}}
\newcommand{\SAm}{\Sigma_{u}^{-}A_i}
\newcommand{\SAp}{\Sigma_{u}^{+}A_i}
 
\newcommand{\soc}[1]{\left\lVert #1 \right\rVert_2} 

\newcommand{\Var}[1]{\text{var} \left[ #1 \right]}

\newcommand{\z}{z_i}

\def\OPFABSW{\mbox{\small{\textsc{\textsf{OPF-SW-AB}}}}}
\def\OPFABNN{\mbox{\small{\textsc{\textsf{OPF-N2N-AB}}}}}

\def\OPFSBSW{\mbox{\small{\textsc{\textsf{OPF-SW-SB}}}}}
\def\OPFSBNN{\mbox{\small{\textsc{\textsf{OPF-N2N-SB}}}}}

\def\GPMABNN{\mbox{\small{\textsc{\textsf{GPM$_i$-N2N-AB}}}}}

\DeclareSIUnit\h{h}
\DeclareSIUnit\A{A}
\DeclareSIUnit\Ah{Ah}
\DeclareSIUnit\MWh{MWh}
\DeclareSIUnit\kWh{kWh}
\DeclareSIUnit\MW{MW}
\DeclareSIUnit\MVAr{MVAr}
\DeclareSIUnit\MVA{MVA}
\DeclareSIUnit\pu{p.u.}
\DeclareSIUnit\EurokWh{\text{\euro}/kWh}
\DeclareSIUnit\USDkWh{\text{\$}/kWh}
\DeclareSIUnit\USDMW{\text{\$}/MW}
\DeclareSIUnit\USDkWhpar{\text{\$}/{kWh^*}}
\DeclareSIUnit\USDMWpar{\text{\$}/{MW^*}}
\DeclareSIUnit\EuroMWh{\text{\euro}/MWh}

\newtheorem{theorem}{\textbf{\textit{Theorem}}}
\theoremstyle{remark}

\newtheorem{proposition}{\textbf{\textit{Proposition}}}

\makeatletter
\renewenvironment{proof}[1][\textbf{\textit{Proof}}]{%
   \par\pushQED{\qed}\normalfont%
   \topsep6\p@\@plus6\p@\relax
   \trivlist\item[\hskip\labelsep\bfseries#1\@addpunct{.}]%
   \ignorespaces
}{%
   \popQED\endtrivlist\@endpefalse
}

\renewcommand\Expected[1]{\mathbb{E}\left[ #1\right]}

\usepackage{mleftright}
\mleftright

\let\originalleft\left
\let\originalright\right
\renewcommand{\left}{\mathopen{}\mathclose\bgroup\originalleft}
\renewcommand{\right}{\aftergroup\egroup\originalright}

\newcolumntype{P}[1]{>{\centering\arraybackslash}p{#1}}

\usepackage{threeparttable}

\usepackage{cellspace} 
\usepackage{cite}

\usepackage{algorithmic}
\usepackage{graphicx}
\usepackage{textcomp}
\usepackage{xcolor}
\def\BibTeX{{\rm B\kern-.05em{\sc i\kern-.025em b}\kern-.08em
    T\kern-.1667em\lower.7ex\hbox{E}\kern-.125emX}}

\begin{document}

\title{Electricity and Reserve Pricing in Chance-Constrained Electricity Markets with Asymmetric Balancing Reserve Policies
}

\author{\IEEEauthorblockN{Alvaro Gonzalez-Castellanos, Anton Hinneck, Robert Mieth, David Pozo, and Yury Dvorkin \vspace{-5mm}}
}

\markboth{Submitted to IEEE Transactions on Power Systems, June 2021}{Author \MakeLowercase{\text
it{et al.}}: Reserve Pricing under Uncertainty}

\maketitle

\begin{abstract}
Recently, chance-constrained stochastic electricity market designs have been proposed to address shortcomings of scenario-based stochastic market designs.
In particular, the use of chance-constrained market-clearing  avoids trading off in-expectation and per-scenario characteristics and yields unique energy and reserves prices.
However, current formulations rely on symmetric control policies based on the aggregated system imbalance, which restricts balancing reserve providers in their energy and reserve commitments. 
This paper extends existing chance-constrained market-clearing formulations by leveraging node-to-node and asymmetric balancing reserve policies and deriving the resulting energy and reserve prices. 
The proposed node-to-node policy allows for relating the remuneration of balancing reserve providers and payment of uncertain resources using a  marginal cost-based approach. Further, we introduce asymmetric balancing reserve policies into the chance-constrained electricity market design and show how this additional degree of freedom affects market outcomes. 
\end{abstract}


\newcommand{\ModelOne}{$\text{CC-POPF-SR}_{\text{I}}$}
\newcommand{\ModelTwo}{$\text{CC-POPF-SR}_{\text{II}}$}
\newcommand{\ModelThree}{$\text{CC-POPF-SR}_{\text{III}}$}

%
\IEEEpeerreviewmaketitle

\section{Introduction}
The deployment of renewable energy sources (RES) challenges the efficiency of wholesale electricity markets, which largely treat RES injections as deterministic and do not internalize their stochasticity rigorously. Although, as Hobbs and Oren discuss in \cite{8606515}, recent market design improvements have targeted (and often succeeded in) improving economic and energy efficiency in the presence of RES stochasticity, they have been “\textit{primarily incremental in nature},” benefiting from enhanced computational capabilities and supply/demand technologies. As a result of these incremental changes, market-clearing procedures have become increasingly complex and market outcomes “\textit{are not transparent and perhaps have contributed to decreases in trading activity}” \cite{8606515}. This lack of transparency inhibits meaningful interpretations of energy and reserve allocations and prices, i.e., there is no technically and economically sound intuition on (i) which resources drive the demand for balancing reserve and how much they should pay for it, and (ii) which resources are most efficient to mitigate this stochasticity and how much they should be paid.  This paper aims to develop a stochastic market design that allows for such interpretations and insights into the energy and reserve price formation under RES stochasticity.

Existing stochastic market designs rely on either scenario-based stochastic programming \cite{papavasiliou2011reserve} or chance-constrained \cite{Bienstock2014} dispatch models, which outperform deterministic benchmarks in terms of the total operating cost and the accuracy of reserve allocations \cite{lubin2016robust}. In \cite{wong2007pricing, pritchard2010single,morales2012pricing}, a two-stage scenario-based stochastic programming framework is used for the day-ahead market-clearing optimization, which yields scenario-specific locational marginal prices (LMPs). Although these LMPs are useful to understand dispatch and price implications of each scenario, it is impossible to ensure cost recovery  (i.e., each producer recovers its production cost from market outcomes) and revenue adequacy (i.e., the payment collected by the market from consumers is greater than the payment by the market to producers) simultaneously in each scenario and in expectation over all scenarios without welfare losses and relying on out-of-market corrections and uplift payments  \cite{kazempour2018stochastic}. 

As an alternative to  \cite{wong2007pricing,pritchard2010single,morales2012pricing,kazempour2018stochastic}, the work in \cite{Kuang2018, Dvorkin2019,Mieth2019,Mieth2020,ratha2019exploring} developed a stochastic market design that internalizes the RES stochasticity by means of its statistical moments (e.g., mean and variance) and chance constraints from \cite{Bienstock2014}. Although stochastic by design, the models in \cite{Kuang2018,Dvorkin2019,Mieth2019,Mieth2020,ratha2019exploring}  yield deterministic reformulations that computationally outperform scenario-based formulations, see \cite{Bienstock2014}, and produce uncertainty- and risk-aware LMPs and reserve prices. Although these prices capture all uncertainty realizations assumed, they are scenario-agnostic, which guarantees cost recovery and revenue adequacy for convex markets \cite{Dvorkin2019,ratha2019exploring}, as well as minimizes the uplift for non-convex markets \cite{Dvorkin2019}. Further, \cite{Dvorkin2019} shows the chance-constrained framework makes it possible to ensure the cost recovery for each uncertainty realization and in expectation without welfare losses. The qualitative analyses in \cite{Dvorkin2019,Mieth2019,Mieth2020} show that these LMPs do not explicitly depend on statistical moments and risk preferences of the market, while reserve prices explicitly depend on these parameters. Despite these computational and market design advantages relative to \cite{wong2007pricing,pritchard2010single,morales2012pricing,kazempour2018stochastic}, the previously developed chance-constrained frameworks have several limitations. First, they typically assume that the RES stochasticity is symmetric, which does not hold in practice, where upward and downward reserve needs in fact vary   \cite{Dvorkin2016, Haupt2019}. Second, while allowing for a nodal reserve allocation, they lack a nodal reserve pricing mechanism, thus preventing from fairly charging and remunerating those resources that drive the need for and provide balancing reserve, respectively.

This paper extends the market design originally proposed in \cite{Kuang2018,Dvorkin2019,Mieth2019,Mieth2020} to accommodate asymmetric reserve provision, node-to-node reserve pricing, and provide techno-economic insights on the energy and reserve price formation process under uncertainty.  Considering the asymmetric reserve provision leads to appropriately sized and allocated reserve requirements drawn from empirical RES statistics (e.g., moments), while the node-to-node reserve pricing mechanism leads to the transparent allocation of (i) uncertainty costs among RES resources and (ii) reserve payments among producers, thus incentivizing the efficient energy and reserve co-optimization to firm up RES generation as necessitated in \cite{8606515}.

\section{Asymmetric Chance-Constrained OPF\label{Sec: Asymmetric OPF}}

\begin{figure}[b]
    \subfloat[Germany, {\footnotesize \texttt{10YDE-EON-----1}}]{
        \includegraphics[width = .46\linewidth]{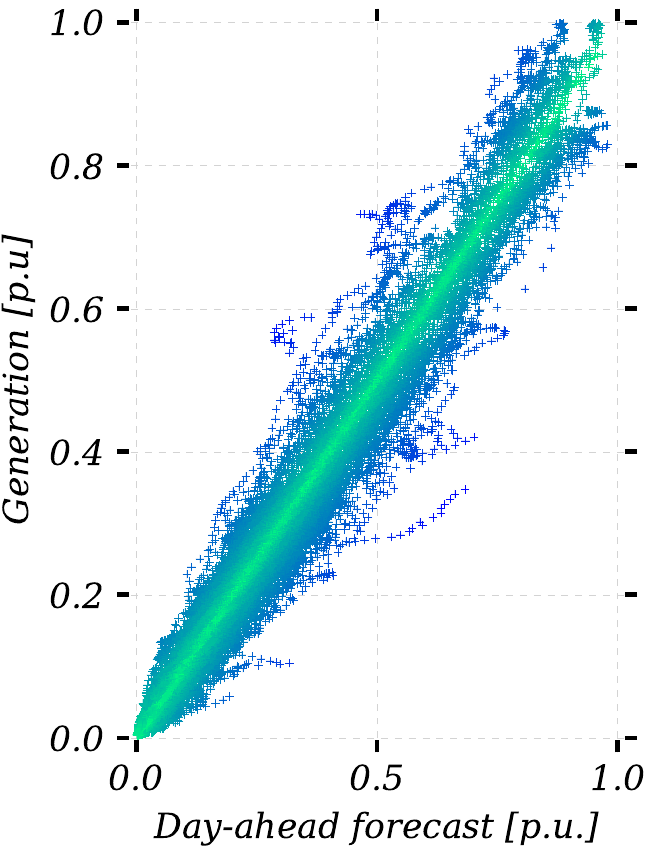}\label{fig:forecasts_ger}}
    \subfloat[Italy, {\footnotesize \texttt{10Y1001A1001A70O}}]{
        \includegraphics[width = .46\linewidth]{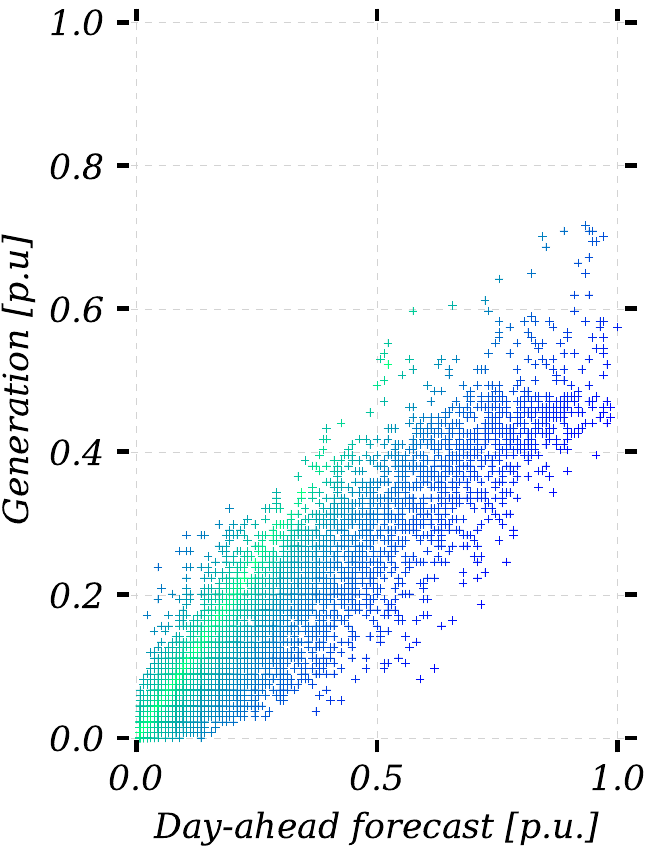}\label{fig:forecasts_it}
    }
  \caption{Day-ahead forecast $w_{t}$ and actual generation $p_{t}$ per unit of maximal power generation in a European aggregation area. An accurate forecast falls on the line $p_{t} = w_{t}$. Figure \ref{fig:forecasts_it} illustrates that forecasting errors in practice can be non-zero-mean.}
\label{fig:forecasts}
\end{figure}

Consider an electricity market operator that uses an optimal power flow (OPF) formulation to compute a least-cost generator schedule. 
As shown in \cite{Kuang2018,Dvorkin2019,Mieth2019,Mieth2020,ratha2019exploring}, traditional deterministic OPF-based market designs (see \cite{Kirschen2018} for reference) can be efficiently robustified against uncertain injections from RES by means of chance constraints. 
For this purpose, injections at every network node $u$ that hosts an uncertain RES resource is modeled as a random variable  ($\textbf{{w}}_u$)  using a given forecast value ($\text{w}_u$) and a random forecast error term ($\om_u$) as follows:
\begin{IEEEeqnarray}{rCll}
    \textbf{{w}}_u &{=}& \text{w}_u{+}\boldsymbol{\omega}_u.
\end{IEEEeqnarray}
Typically, $\om_u$ is assumed to be zero-mean and normally distributed, \cite{Bienstock2014,lubin2016robust,Dvorkin2020,Mieth2019,Mieth2020,Kuang2018, ratha2019exploring}. 
However, in practice, empirical measurements of solar and wind power forecast errors are often asymmetric and can not be captured well by a  normal distribution, \cite{Dvorkin2016, Haupt2019}. 
As a result, the impact of forecast errors from nodes with significant asymmetries might be over- or under-estimated. 
For example, by inspecting day-ahead forecast and actual generation data recorded at two nodes in the ENTSO-E (European Network of Transmission System Operators for Electricity) in 2020, plotted in Figure~\ref{fig:forecasts}, various deviations of the forecast error from a normal distribution become evident. 
While the forecast at the German node, Figure~\ref{fig:forecasts}(a), is symmetric and may be parametrized as a normal distribution, the forecast at the Italian node, Figure~\ref{fig:forecasts}(b), is strongly skewed to overestimate generator output and a normal distribution is unsuitable to model the forecast error distribution.
Additionally, there exists a tendency of RES operators to report generation values lower than the day-ahead forecast; which could be a result of strategic forecast offering by the generator in order to obtain additional benefits from complex and not always transparent market regulations. See in \cite[Section 6.6.]{conejo2010decision}.

As a result, assuming normally distributed (symmetric) forecast errors in combination with a symmetric balancing reserve policy as in \cite{Bienstock2014,lubin2016robust,roald2013analytical,xie2017distributionally,Dvorkin2020,Mieth2019,Mieth2020,Roald2016,Fang2019,Kuang2018,chertkov2017chance,Bienstock2018} can lead to ineffective and inefficient operating decisions and electricity prices. 
Therefore, to adequately capture possible forecast error asymmetries and improve the efficiency of balancing reserve quantification and allocation, this paper explicitly models negative and positive forecast errors (i.e., real-time energy deficit and surplus, respectively) and proposes an asymmetric, node-to-node balancing reserve policy. 
Sections~\ref{ssec:UncertaintyModelling} and \ref{ssec:asymmetric_ccopf} below describe the asymmetric model of uncertain RES injections and introduce the proposed balancing reserve policy, respectively. Section~\ref{ssec:model_formulation} derives the resulting market-clearing optimization. 

\subsection{Asymmetric Uncertainty Model}\label{ssec:UncertaintyModelling}

\begin{figure}[t]
\centering
    \subfloat[Symmetric forecast error]
        {\includegraphics[width=0.45\columnwidth]{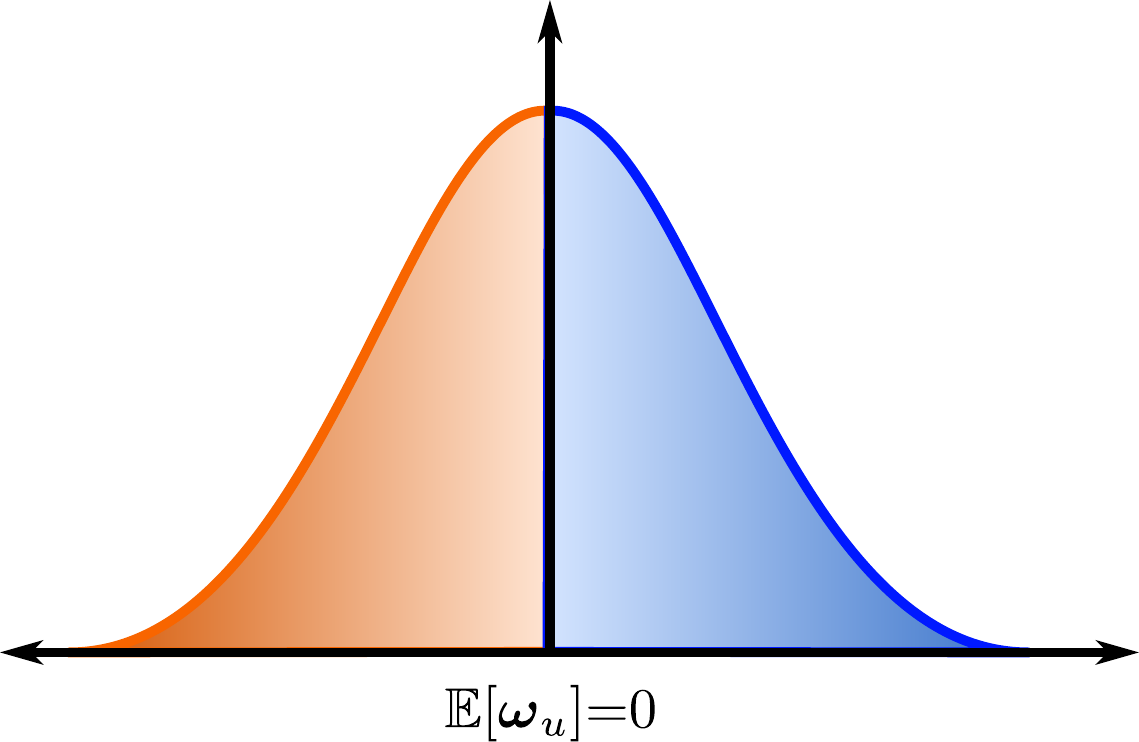}
  		\label{Fig: SymmetricDist}}
    \subfloat[Asymmetric forecast error]
        {\includegraphics[width=0.45\columnwidth]{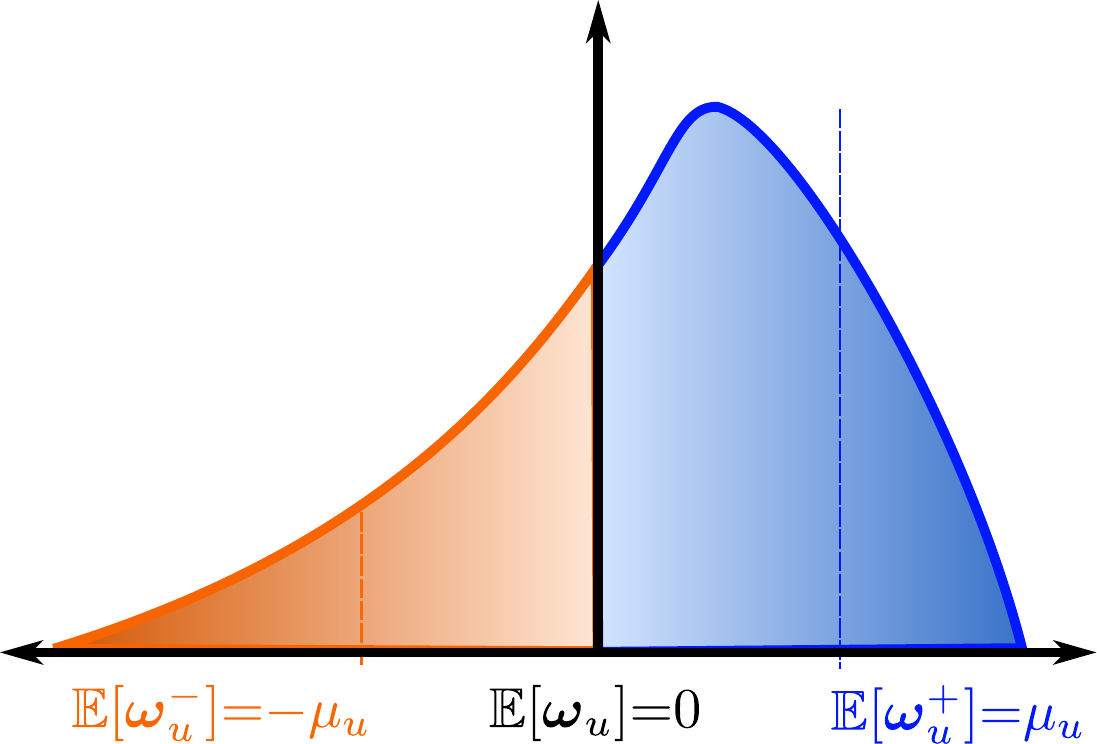} \label{Fig: AsymmetricDist}}
   \caption{Nodal forecast errors \(\om_u\) for  (a) symmetric or (b) asymmetric probability distributions with zero expected value.  \label{Fig: CompDistrib}}
\end{figure}

Consider random variable $\om_u$:
\begin{equation}
 \om_u = \om^\text{-}_u + \om^\text{+}_u,
 \label{eq:split_omu}
\end{equation}
where $\om^\text{-}_u \leq 0 $ and $\om^\text{+}_u \geq 0$ are the negative and positive components of the forecast error such that $\om^\text{-}_u \om^\text{+}_u = 0$, i.e., $\om^\text{-}_u$ and $\om^\text{+}_u$ are mutually exclusive events.
First, we assume that any expected systematic forecast error can be considered as a fixed parameter in the model (e.g., by attributing it to fixed load forecasts), so that the total imbalance ($\om_u$) can always be corrected to be zero-mean. 
Second, from \eqref{eq:split_omu} and the linearity  of the expectation operator, it follows that $\Expected{\om_u} = \Expected{\om_u^- + \om_u^+ } = 0$, i.e., $\Expected{\om^\text{+}_u} = - \Expected{\om^\text{-}_u} = \mu_u$.
Figures~\ref{Fig: CompDistrib}(a)--(b) illustrate this representation, which allows for modeling asymmetric distributions.

In turn, mean $\mu_u$ of the forecast error distributions can be obtained by considering them as distributions truncated at zero. Thus, for an estimated continuous distribution, it follows:
\begin{IEEEeqnarray}{rl}
    \mu_u = \Expected{\om_u {\geq} 0} 
    =
    \Expected{\om^\text{+}_u} = 
    {\int_{0}^{\infty}} {\omega_u f_{\om_u}
    \left(\omega_u\right)} ~d\omega_u, \quad &\forall u \IEEEeqnarraynumspace
\end{IEEEeqnarray}
where $f_{\om_u}$ is the probability density function of $\om_u$. 

\subsection{Asymmetric Node-to-Node Reserve Balancing}
\label{ssec:asymmetric_ccopf}

To ensure that the power system remains balanced (i.e., generation equals demand at all times), controllable generators with balancing reserve capabilities react to any real-time imbalances $\om^\text{-}_u$ or $\om^\text{+}_u$ by increasing or decreasing their power outputs. 
We model this balancing reserve policy using asymmetric ``node-to-node'' balancing participation factors $\alpha_{iu}^- \in [0,1]$ and $\alpha_{iu}^+ \in [0,1]$ that model the participation of a generator at bus $i$ in compensating a negative and positive imbalance caused by a stochastic RES at bus $u$. Note that for ease of notation, we assume that each node $i$ hosts exactly one controllable generator and each node $u$ hosts exactly one uncertain RES.
The resulting uncertain power output $\boldsymbol{p}_{i}$ of each controllable generator is then given by:
\begin{IEEEeqnarray}{lcll}
    \boldsymbol{p}_{i} = p_{i} {-} \sum_u \left( \alpha^\text{-}_{iu} \om_u^\text{-} {+} \alpha^\text{+}_{iu}\om_u^\text{+} \right), \qquad & \forall i \IEEEyesnumber\IEEEyessubnumber*\label{Eq: p AB N2N}\\
    \sum_i\alpha_{iu}^\text{-} = 1, \quad  \sum_i\alpha_{iu}^\text{+} = 1, & \forall u \label{Eq: Sufficiente a_iu AB}
\end{IEEEeqnarray} 
where $p_i$ and $\alpha_{iu}^\text{-}, \alpha_{iu}^\text{+} \in [0,1]$ are decision variables that define the power output  and asymmetric  participation factors for balancing reserve provision. Constraints \eqref{Eq: Sufficiente a_iu AB} ensure system balance for all outcomes $\om^\text{-}_u$, $\om^\text{+}_u$ \cite{Bienstock2018,Roald2016}. Note that random variables are written in \textbf{bold} font throughout this paper.

\subsection{Model Formulation}
\label{ssec:model_formulation}

Following \cite{Kuang2018, Dvorkin2020,Mieth2019,Mieth2020, ratha2019exploring}, the market objective is to minimize the expected production cost of conventional generators.  By replacing $\p_i$ from \eqref{Eq: p AB N2N} in \eqref{Eq: Gen c_i}, the expected generation cost $C_i$ of generator $i$ can be calculated as in Eq.~\eqref{Eq: GCost AB N2N} (see Appendix \ref{App: Asymmetric cost N2N} for details):
%
\begin{IEEEeqnarray}{rl}
    {c_i}{(\boldsymbol{p}_{i})} &= {c_{2,i}} {{\left( \boldsymbol{p}_{i} \right)}^{2}} {+} c_{1,i} {\left( \p_{i} \right)} {+} c_{0,i}, \quad \forall i \label{Eq: Gen c_i} \\
   \mathbb{E}[c_i(\boldsymbol{p}_{i})] & = C_i  
    = {c_{2,i}} p_i^2 + c_{1,i} p_i + c_{0,i}  \: - \label{Eq: GCost AB N2N}\\
    & \hphantom{=}  {\left(2 c_{2,i} p_i  {+} c_{1,i} \right)}{\left( M A_i \right)}{+} c_{2,i} {\left[ \left(M {\cdot} A_i\right)^2 {+} S_i^2 \right]}{,} ~~  \forall i, \nonumber
\end{IEEEeqnarray}
where ${M}$ is the row vector of expected negative and positive forecast errors, ${M=\Expected{\om}}$, with $\boldsymbol{\omega}{=}\left[\om_{1}^-,\om_{2}^-,...,\om_{|\mathcal{U}|}^-,\om_{1}^+,\om_{2}^+,...,\om_{|\mathcal{U}|}^+ \right]$, and $A_i$ is the vector of balancing participation factors for generator $i$, $A_i{=}\left[\alpha_{i1}^-,\alpha_{i2}^-,...,\alpha_{i|\mathcal{U}|}^-, \alpha_{i1}^+, \alpha_{i2}^+,..., \alpha_{i|\mathcal{U}|}^+ \right]^{\top}$. The $\cdot$ operator denotes a dot product of two vectors. The covariance matrix of the forecast errors is given by $\Sigma{=}cov\left[\Om\right]$ and $S_i {=} \sqrt{A_i^{\!\top} \Sigma A_i} {=}\AZ$.
%

Upper and lower power output limits are enforced by the chance constraints \eqref{Eq: PmaxI} and \eqref{Eq: PminI}, respectively, and parameter $\epsilon_i\leq 0.5$ is chosen small enough to ensure constraint satisfaction  with an $1{-}\epsilon_i$ guarantee:
\begin{subequations}
\label{Eq: PminmaxI}
\begin{IEEEeqnarray}{rCll}
    \mathbb{P}\left[\boldsymbol{p}_{i}\leq \overline{P}_{i}\right]  &\geq& 1-\epsilon_{i}, &\forall i\IEEEeqnarraynumspace\label{Eq: PmaxI}\\
    \mathbb{P}\left[-\boldsymbol{p}_{i}\leq -\underline{P}_{i}\right] &\geq& 1-\epsilon_{i}, \quad &\forall i. \IEEEeqnarraynumspace\label{Eq: PminI}
\end{IEEEeqnarray}
\end{subequations}
Constraints \eqref{Eq: PmaxI} and \eqref{Eq: PminI} enforce that the probability of generator $i$ exceeding its maximum and minimum capacity limits ($\overline{P}_i$ and $\underline{P}_i$ ),
while employing a pre-defined balancing reserve policy, is tolerable by the market operator.
Therefore, the underlying probability distribution functions representing the different RES forecast errors must be accounted for accurately.
However, unlike for Gaussian distributions, the aggregation of asymmetric distributions does not have an equivalent deterministic reformulation.
To overcome this issue, we apply the Chebyshev inequality to \eqref{Eq: PminmaxI} for representing its satisfaction probability for any kind of underlying distribution \cite{Dvorkin2020}.
Thus, the capacity constraints \eqref{Eq: PminmaxI} can be reformulated by means of the Chebyshev approximation, as the following equivalent and deterministic constraints in \eqref{Eq: change constraint Max sqrt} and \eqref{Eq: change constraint Min sqrt}:
\begin{IEEEeqnarray}{rCl"l}
    p_{i} 
    - {{M} {\cdot} A_i}
    + \z {S_i}
     &\leq& 
    \overline{P}_{i}, & \forall i \IEEEyesnumber\IEEEyessubnumber\label{Eq: change constraint Max sqrt} \\
       -p_{i} 
    + {{M} {\cdot} A_i}
    + \z {S_i}
    &\leq&   
    -\underline{P}_{i}, & \forall i, \IEEEyessubnumber\label{Eq: change constraint Min sqrt} 
\end{IEEEeqnarray}
where the parameter $\z {=} \sqrt{\left(1-\epsilon_i\right)/\epsilon_i}$.

Using \eqref{Eq: Sufficiente a_iu AB}, \eqref{Eq: GCost AB N2N}, \eqref{Eq: change constraint Max sqrt}, and \eqref{Eq: change constraint Min sqrt}, the OPF problem with a deterministic equivalent of the chance-constrained generator dispatch and node-to-node asymmetric balancing reserve policy (\OPFABNN{}) is formulated in Model \ref{OPFABNN}. The objective function  in  \eqref{Eq: Obj OPFABNN} minimizes the expected system generation cost and is subject to dc power flow constraints \eqref{Eq: P Balance}-\eqref{Eq: Fmin}, asymmetric node-to-node balancing reserve adequacy requirements \eqref{Eq: Alpha sum m u}--\eqref{Eq: Alpha sum p u}, and chance-constrained generator output limits \eqref{Eq: Pmax u AB}--\eqref{Eq: Pmin u AB}.
In the \OPFABNN{}, the network buses are indexed by $i$, the RES buses by $u \in \mathcal{U}_i$, and the power lines by the tuple $\left(i,j\right)$ with $j \in \mathcal{C}_i$; where $\mathcal{C}_i$ is the set of nodes directly connected to node $i$.
The nodal demand is given by $D_i$ and the line reactance by $X_{ij}$.
The line flow is represented by $f_{ij}$ and the nodal voltage angle by $\theta_i$.
Line capacity limits are set to $\overline{F}_{ij}$. 

The decision variables of \OPFABNN{} are collected in set  $\Xi = \{ p_i, \alpha_{iu}^\text{-}, \alpha_{iu}^\text{+}, f_{ij}, \theta_i \}$. Greek symbols, in parenthesis, define dual multipliers associated with each constraint. 
Note that \OPFABNN{} is a non-linear program (NLP) with second-order conic constraints. Specifically, objective \eqref{Eq: Obj OPFABNN} contains bi-linear terms, which are dealt with in Section~\ref{Sec: Case Study}.
However, the second-order conic solution space formed by constraints \eqref{Eq: P Balance}--\eqref{Eq: Pmin u AB} is convex and, thus, enables the analytical derivations in the following sections. 

%
\begin{model}[t]
\caption{\OPFABNN }
\label{OPFABNN}
\begin{subequations} 
\label{Mod: OPFABNN}
\vspace{-2\jot}
\begin{IEEEeqnarray}{lll}
    \minimize_\Xi &~
    \sum_i \! \Bigg[&
    {c_{2,i}} {{ {p}_{i}}^{2}} + c_{1,i} {\left( p_{i} \right)} + c_{0,i} \IEEEnonumber\vspace{-2\jot}\\
    &&~{-} {\left(2 c_{2,i} p_i{+} c_{1,i}\right)}
    {\left(M \cdot A_i\right)} {+} c_{2,i} {\!\left[\left(M \cdot A_i\right)^2 {+} S_i^2 \right]\Bigg]} \IEEEnonumber\\
    \label{Eq: Obj OPFABNN}
    \vspace{-2\jot}
    \\
    \IEEEeqnarraymulticol{3}{l}{\text{subject to:}} \IEEEnonumber \vspace{-2\jot}
\end{IEEEeqnarray}
\begin{IEEEeqnarray}{rlrrrl}
    (\lambda^\text{}_i){:} \quad & p_{i} + 
    \!\textstyle{\sum_{u \in \mathcal{U}_i}} \text{w}_u - D_{i}^\text{} - \!\textstyle{\sum_{j \in \mathcal{C}_i}} f_{ij}^{} = 0, &&&\forall i \IEEEeqnarraynumspace \label{Eq: P Balance} \\
    (\psi^\text{}_{ij}){:}  \quad& f_{ij} - \frac{1}{X_{ij}} \left(\theta_i - \theta_j\right) = 0,\quad & \IEEEeqnarraymulticol{3}{r}{\forall i, j {\in} \mathcal{C}_i} \IEEEeqnarraynumspace\label{Eq: F def}\\
    (\overline{\eta}_{ij}){:} \quad & f_{ij}\leq \overline{F}_{ij}, &\IEEEeqnarraymulticol{3}{r}{\forall i, j {\in} \mathcal{C}_i} \IEEEeqnarraynumspace \label{Eq: Fmax}\\
    (\underline{\eta}_{ij}){:} \quad & \hspace{-2.5\jot}-f_{ij}\leq -\overline{F}_{ij}, &\IEEEeqnarraymulticol{3}{r}{\forall i, j {\in} \mathcal{C}_i} \IEEEeqnarraynumspace\label{Eq: Fmin}\\
    (\chi_u^\text{-}){:} \quad &  \!\sum_i\alpha_{iu}^\text{-} = 1, &&&\forall u \IEEEeqnarraynumspace \label{Eq: Alpha sum m u}\\
    (\chi_u^\text{+}){:} \quad &  \!\sum_i\alpha_{iu}^\text{+} = 1, &&&\forall u \IEEEeqnarraynumspace \label{Eq: Alpha sum p u}\\
    (\overline{\delta}_{i}){:} \quad 
    & \IEEEeqnarraymulticol{2}{l}{p_{i} 
    {-} {{M} {\cdot} A_i}
    {+} \z{S_i}
    {\leq} 
    \hspace{2.5\jot}\overline{P}_{i},} &&\forall i \IEEEeqnarraynumspace \label{Eq: Pmax u AB}\\
    (\underline{\delta}_{i}){:} \quad 
    & \IEEEeqnarraymulticol{2}{l}{\hspace{-2.5\jot}-p_{i} 
    {+} {{M} {\cdot} A_i}
    {+} \z{S_i}
    {\leq} 
    {-}\underline{P}_{i},} &&\forall i\IEEEeqnarraynumspace &\hspace{-5.5\jot}. \label{Eq: Pmin u AB}
    \vspace{-2\jot}
\end{IEEEeqnarray}
\end{subequations}
\end{model}

\section{Pricing Energy and Balancing Reserve Provision}
\label{Sec: Prices AB NN}
From the perspective of an electricity market operator, \OPFABNN{} resembles a market-clearing problem computing the optimal (least-cost) generator dispatch and balancing reserve decisions $p_i^\ast$, $\alpha^{\text{-},\ast}_{iu}$, $\alpha^{\text{+},\ast}_{iu}$ and efficient prices for energy $\pi^\text{p}$ and reserves $\pi_u^\text{-}$,  $\pi_u^\text{+}$.
More specifically, the price-quantity tuples $(\pi^\text{p}, p_i^\ast)$, $(\pi_u^\text{-}, \alpha^{\text{-}})$, $(\pi_u^\text{+}, \alpha^{\text{+},\ast}_{iu})$ are efficient, if they support a competitive equilibrium, i.e., they ensure \cite{Dvorkin2020}:
\begin{enumerate}
    \item The market clears at $p_{i}^* + \!\textstyle{\sum_{u \in \mathcal{U}}} \text{w}_u  - \!\sum_{j \in \mathcal{C}_i} f_{ij}^{} {=} D_{i},~ \forall i$, $\sum_i \alpha_{iu}^{\text{-},*} {=} \sum_i \alpha_{iu}^{\text{+},*}=1,~ \forall u$.
    
    \item  Prices $\pi^\text{p}$, $\pi_u^\text{-}$ and quantities $p_i^\ast$, $\alpha^{\text{-},\ast}_{iu}$ maximize the profit of $i$ given by $\Pi_i = \pi_{i}^\text{p} p_{i}^* {+} {\sum_{u}{\!\left(\pi^\text{-}_{u}\alpha_{iu}^{\text{-},*} {+} \pi^\text{+}_{u}\alpha_{iu}^{\text{+},*} \right)}} - C_i$, i.e., there is no incentive to deviate these market outcomes.
\end{enumerate}

To derive these prices from \OPFABNN{}, we use duality theory similar to the previous work on stochastic market designs with chance constraints \cite{Dvorkin2020,Mieth2019,Mieth2020}. However, unlike in \cite{Dvorkin2020,Mieth2019,Mieth2020}, the prices derived from \OPFABNN{} will support asymmetric balancing reserve needs and  account for node-to-node participation factors, thus enabling a more precise allocation and deployment of balancing reserve. 
Thus, we begin with the following proposition:
%
\begin{proposition}\label{Prop: AB N2N}
    Consider the \OPFABNN{} model. Let $\lambda_i$, $\chi_u^\text{-}$ $\chi_u^\text{+}$ be the locational marginal prices (LMPs) for energy and asymmetric locational balancing prices (A-LBPs) defined as dual multipliers of constraints \eqref{Eq: P Balance}, \eqref{Eq: Alpha sum m u}, and \eqref{Eq: Alpha sum p u}, respectively. Then, these prices can be expressed as:
    \begin{IEEEeqnarray}{lll}
        \IEEEeqnarraymulticol{2}{l}{\lambda_i {=}  2c_{2,i}\left(p_i-M{\cdot}A_i\right) + c_{1,i} + \overline{\delta}_{i} - \underline{\delta}_{i}, }& {\forall i} \IEEEeqnarraynumspace \IEEEyesnumber\IEEEyessubnumber*\label{Eq: KKT lambda_i AB NN}\\
        %
        \IEEEeqnarraymulticol{2}{l}{
        \chi_u^\text{-} {=} \frac{1}{|\mathcal{G}|} \!\sum_{i} \!\left[{ \mu_u \lambda_i {+} \left(\SAm\right)\left(2c_{2,i}{+}\frac{\z}{S_i}\left(\overline{\delta}_i {+}\underline{\delta}_i\right)\right) }\right]\!{,}~}&\forall u \IEEEeqnarraynumspace\label{Eq: KKT chi_u AB NN m}\\
        %
        \IEEEeqnarraymulticol{2}{l}{
        \chi_u^\text{+} {=} \frac{1}{|\mathcal{G}|} \!\sum_{i} \!\left[{ -\mu_u \lambda_i {+} \left(\SAp\right)\left(2c_{2,i}{+}\frac{\z}{S_i}\left(\overline{\delta}_i {+}\underline{\delta}_i\right)\right) }\right]\!{,}~}&\forall u{.} \IEEEeqnarraynumspace\label{Eq: KKT chi_u AB NN p}
    \end{IEEEeqnarray}
where $\Sigma_{u}^{-}$ and $\Sigma_{u}^{+}$ are the rows of the covariance matrix $\Sigma$ related to the asymmetric forecast errors $\om_u^\text{-}$ and $\om_u^\text{+}$, respectively.\footnote{For example, in a system with two RES ${\Sigma_{1}^{-} = \left[ Var(\om_1^\text{-}), cov(\om_1^\text{-},\om_1^\text{+}), cov(\om_1^\text{-},\om_2^\text{-}), cov(\om_1^\text{-},\om_2^\text{+})\right]}$.}
\end{proposition}
\begin{proof}
Consider the first-order optimality conditions of \eqref{Mod: OPFABNN} for variables $p_i$, $\alpha_{iu}^\text{-}$, and $\alpha_{iu}^\text{+}$:
\begin{IEEEeqnarray}{lr}
    2c_{2,i}\left(p_i-M{\cdot}A_i\right) {+} \overline{\delta}_{i} {-} \underline{\delta}_{i}{-}\lambda_i = 0,~&{\forall i}
    \IEEEyesnumber \IEEEyessubnumber*\IEEEeqnarraynumspace \label{Eq: KKT p_i AB NN}\\
    \IEEEeqnarraymulticol{2}{l}{
    \left(\SAm\right)\left(2c_{2,i}{+}\frac{\z}{S_i}\left(\overline{\delta}_i {+}\underline{\delta}_i\right)\right) 
    {+} \mu_u \left(2c_{2,i}\left(p_i{-}M{\cdot}A_i\right) {+} \overline{\delta}_{i} {-} \underline{\delta}_{i}\right) } \IEEEnonumber\\
    &{-}\chi_u^\text{-}=0, {\forall i{,}u}  \IEEEeqnarraynumspace
    \label{Eq: KKT alpha_iu AB NN m}\\
    \IEEEeqnarraymulticol{2}{l}{
    \left(\SAm\right)\left(2c_{2,i}{+}\frac{\z}{S_i}\left(\overline{\delta}_i {+}\underline{\delta}_i\right)\right) 
    {-}\mu_u \left(2c_{2,i}\left(p_i{-}M{\cdot}A_i\right) {+} \overline{\delta}_{i} {-} \underline{\delta}_{i}\right)} \IEEEnonumber\\
    &{-}\chi_u^\text{+}=0, {\forall i{,}u.}  \IEEEeqnarraynumspace
    \label{Eq: KKT alpha_iu AB NN p}
\end{IEEEeqnarray}
Expression \eqref{Eq: KKT lambda_i AB NN} immediately follows from separating $\lambda_i$ from \eqref{Eq: KKT p_i AB NN}. Note that $\chi_u^\text{-}$ and  $\chi_u^\text{+}$ are the dual variables of the nodal balancing reserve adequacy constraints in Eq.~\eqref{Eq: Alpha sum m u} and \eqref{Eq: Alpha sum p u} and, therefore, are specific to the balancing needs of node \(u\). 
Note that expressions \eqref{Eq: KKT alpha_iu AB NN m} and \eqref{Eq: KKT alpha_iu AB NN p} have the same value for each generator $i$ and can be summed over  $|\mathcal{G}|$ generators, leading to $\sum_{i \in \mathcal G} \chi_u^\text{-} = |\mathcal{G}| \chi_u^\text{-} $ and  $\sum_{i \in \mathcal G} \chi_u^\text{+} = |\mathcal{G}| \chi_u^\text{+}$. Hence,, dividing \eqref{Eq: KKT alpha_iu AB NN m} and \eqref{Eq: KKT alpha_iu AB NN m}  by   
\(|\mathcal{G}|\) leads to  \eqref{Eq: KKT chi_u AB NN m} and \eqref{Eq: KKT chi_u AB NN p}.
\qedhere
\end{proof}
We now prove that the dual variables used as prices in Proposition~\ref{Prop: AB N2N} and the primal dispatch decisions satisfy the conditions for a competitive equilibrium.

\begin{theorem}\label{Th: AB N2N}
\textup{Consider the \OPFABNN{} model and let $\lambda_i^*$, $\chi^{-,*}_u$ and $\chi^{+,*}_u$ denote the optimal value of the dual multipliers. Using  $\pi^\text{p}=\lambda_i^*$, $\pi_{u}^\text{-}=\chi^{-,*}_u$ and $\pi_{u}^\text{+}=\chi^{+,*}_u$ as the LMPs and A-LBPs constitutes a competitive energy and balancing reserve equilibrium prices.
}
\end{theorem}
\begin{proof}

Each generator $i$ is modeled using a profit-maximizing problem (GPM$_i$) given in Model~\ref{Mod.GPMABNN}, which is abbreviated below as \GPMABNN. 
\begin{model}[t]
\caption{\GPMABNN}
\label{Mod.GPMABNN}
\vspace{-5pt}
\begin{IEEEeqnarray}{lllCll}
    & \IEEEeqnarraymulticol{4}{l}{\maximize \quad\Pi_i = \pi_{i}^\text{p} p_{i} {+} {\!{\sum_{u}{\!\left(\pi^\text{-}_{u}\alpha_{iu}^\text{-} {+} \pi^\text{+}_{u}\alpha_{iu}^\text{+} \right)}}} {-} C_i\quad} \IEEEeqnarraynumspace  
    \IEEEyesnumber\label{Eq: Profit max AB NN}\IEEEyessubnumber*\\
    & \text{subject to:}  \IEEEnonumber\\
    \left(\overline{\delta}_{i}\right){:} \quad 
    & p_{i} 
    {-} {M {\cdot} A_i}
    {+} \z{S_i}
    \leq 
    \hspace{2.5\jot}\overline{P}_{i},\qquad &\forall i \IEEEeqnarraynumspace \label{Eq: Pmax gen u AB}\\
    \left(\underline{\delta}_{i}\right){:} \quad 
    & \hspace{-2.5\jot}-p_{i} 
    {+} {M {\cdot} A_i}
    {+} \z{S_i}
    \leq  
    {-}\underline{P}_{i},\qquad &\forall i.\IEEEeqnarraynumspace \label{Eq: Pmin gen u AB}
\end{IEEEeqnarray}
\end{model}
The first-order optimality conditions of this optimization with respect to variables $p_i$, $\alpha_{iu}^\text{-}$, and $\alpha_{iu}^\text{+}$ lead to the expressions \eqref{Eq: KKT p_i AB NN price},\eqref{Eq: KKT alpha_iu AB NN m price}, and \eqref{Eq: KKT alpha_iu AB NN p price}, respectively:
\begin{IEEEeqnarray}{lr}
    {\pi_i =  2c_{2,i}\left(p_i-M{\cdot}A_i\right) + c_{1,i} + \overline{\delta}_{i} - \underline{\delta}_{i},~}\hspace{1.5cm}&{\forall i}
    \IEEEyesnumber \IEEEyessubnumber*\IEEEeqnarraynumspace \label{Eq: KKT p_i AB NN price}\\
    \IEEEeqnarraymulticol{2}{l}{\pi_u^\text{-} =
    {\mu_u \left(2c_{2,i}\left(p_i{-}M{\cdot}A_i\right) {+} c_{1,i} {+} \overline{\delta}_{i} {-} \underline{\delta}_{i}\right) }}\IEEEnonumber\\
    \IEEEeqnarraymulticol{2}{r}{{+} \left(\SAm\right)\left(2c_{2,i}{+}\frac{\z}{S_i}\left(\overline{\delta}_i {+}\underline{\delta}_i\right)\right),\quad {\forall i{,}u.}} \IEEEeqnarraynumspace\label{Eq: KKT alpha_iu AB NN m price}\vspace{2\jot}\\
    \IEEEeqnarraymulticol{2}{l}{\pi_u^\text{+} =
    {{-}\mu_u \left(2c_{2,i}\left(p_i{-}M{\cdot}A_i\right) {+} c_{1,i} {+}  \overline{\delta}_{i} {-} \underline{\delta}_{i}\right)  },}\IEEEnonumber\\
    \IEEEeqnarraymulticol{2}{r}{{+} \left(\SAp\right)\left(2c_{2,i}{+}\frac{\z}{S_i}\left(\overline{\delta}_i {+}\underline{\delta}_i\right)\right),\quad {\forall i{,}u.}} \IEEEeqnarraynumspace\label{Eq: KKT alpha_iu AB NN p price}
\end{IEEEeqnarray}
Due to the equality of first-order optimality conditions in \eqref{Eq: KKT p_i AB NN price}--\eqref{Eq: KKT alpha_iu AB NN p price} and in \eqref{Eq: KKT lambda_i AB NN}--\eqref{Eq: KKT chi_u AB NN p}, it follows that the primal ($p_i, \alpha^\text{+}_{iu}, \alpha^\text{-}_{iu}$) and dual ($\lambda_i, \chi_u^\text{-}, \chi_u^\text{+}$) solution from the \OPFABNN{} model, which satisfies constraints \eqref{Eq: P Balance}, \eqref{Eq: Alpha sum p u} and \eqref{Eq: Alpha sum m u}, also solves the \GPMABNN{} model for each generator, which maximizes $\Pi_i$. Hence, the prices obtained in Proposition \ref{Prop: AB N2N} and derived from the \OPFABNN, \eqref{Eq: KKT lambda_i AB NN}--\eqref{Eq: KKT chi_u AB NN p}, constitute a competitive equilibrium, i.e. $\pi^\text{p}=\lambda_i^*$, $\pi_{u}^\text{-}=\chi^{-,*}_u$ and $\pi_{u}^\text{+}=\chi^{+,*}_u$.
\qedhere
\end{proof}
\section{Stochastic Markets under Symmetric and System-wide Balancing Reserve Policies} \label{Sec:modifications}
Using the results derived in Section~\ref{Sec: Prices AB NN}, we consider several particular cases that allow for relating the obtained prices to existing results in \cite{lubin2016robust, Bienstock2014, Kuang2018, Mieth2019, Mieth2020, Dvorkin2020,Bienstock2018,Roald2016 }. 
A summary of the relationship between the pricing for the different balancing reserve policies is provided in Table \ref{Tab: Gen values}.

\subsection{Asymmetric System-wide Balancing -- \OPFABSW{} \label{Sec: ASB}}
Instead of the node-to-node balancing reserve policy \eqref{Eq: p AB N2N}, the proposed \OPFABNN{} can be modified to use a balancing reserve policy based on the aggregated system-wide uncertainty as in, e.g., \cite{Kuang2018, Mieth2019, Mieth2020, Dvorkin2020}.
Thus, the system-wide uncertainty is given as:
\begin{IEEEeqnarray}{C}
    \Om^\text{-}{=}\sum_u{\boldsymbol{\omega}^\text{-}_u}, \quad 
    \Om^\text{+}{=}\sum_u{\boldsymbol{\omega}^\text{+}_u} \label{Eq: expected Om_pm},\\
    \noalign{\noindent where in turn $\Om^\text{+}$ and $\Om^\text{-}$ are related as:
    \vspace{2\jot}}
    \Expected{\Om^\text{+} {+} \Om^\text{-} } {=} \Expected{\Om} {=} 0, \label{Eq: expected Om}
\end{IEEEeqnarray}
and the the node-to-node participation factors used in the \OPFABNN{} are set to system-wide participation factors by enforcing $\alpha_{iu}^\text{-} {=} \alpha_i^\text{-}$ and $\alpha_{iu}^\text{+} {=} \alpha_i^\text{+}$, which leads to $A_i {=} \left[\alpha_i^\text{-}, \alpha_i^\text{+} \right]^\T$.
Although this approach resembles the pricing mechanism from \cite{Kuang2018, Mieth2019, Mieth2020, Dvorkin2020}, there is a notable distinction of this representation, because it still accounts for asymmetric distributions and asymmetric provision of balancing reserve.

We define $\Om^\text{-/+} = \left[\Om^\text{-}, \Om^\text{+}\right]$ and its expected value as $M = \Expected{\Om^\text{-/+}} = \left[-\nu, \nu\right]$, where $\nu = \sum_u \mu_u = \sum_u \Expected{\om_u^\text{+}}= {-}\sum_u \Expected{\om_u^\text{-}}$ is the expected value of the aggregated asymmetric forecast errors.
Finally, the covariance matrix of $\Om^\text{-/+}$ is defined as $\Sigma=cov\left(\Om^\text{-/+}\right)$.

\begin{IEEEeqnarray}{lrCll}
    \noalign{\textit{1) Generation Output:}
    Under the asymmetric balancing reserve policy with system-wide participation factors, we replace  Eqs.~\eqref{Eq: p AB N2N} and \eqref{Eq: Sufficiente a_iu AB} with: \vspace{-4\jot}}
    \IEEEnonumber\\
    &\IEEEeqnarraymulticol{3}{l}{\boldsymbol{p}_{i} = p_{i} {-} \alpha^\text{-}_i \mathbf{\Omega^\text{-}} {-} \alpha^\text{+}_i \mathbf{\Omega^\text{+}},} & \forall i \IEEEyesnumber\IEEEyessubnumber* \label{Eq: p AB}\\
    \left(\chi^\text{-},\chi^\text{+}\right){:}\quad &\IEEEeqnarraymulticol{3}{l}{\hspace{0.5\jot} 1 \hspace{\jot} = \sum_i \alpha_i^\text{-} = \sum_i \alpha_i^\text{+}.}  \\
    \noalign{\vspace{2\jot}\textit{2) Energy and Balancing Prices:}
    Following the procedure described in Proposition~\ref{Prop: AB N2N}, we derive the following prices: 
    \vspace{-3\jot}} \IEEEnonumber\\
    \IEEEeqnarraymulticol{4}{l}{
    \lambda_i = c_{2,i}\left(p_i-M{\cdot}A_i\right) {+} c_{1,i} {+} \overline{\delta}_{i} {-} \underline{\delta}_{i},}&\forall i \IEEEeqnarraynumspace\label{Eq: KKT lambda_i AB SW}\\
    \IEEEeqnarraymulticol{4}{l}{\chi^\text{-} \hspace{0.75\jot} = \frac{1}{|\mathcal{G}|} \!\sum_{i} \!\left[{ \nu \lambda_i {+} \left(\Sigma^\text{-}A_i\right)\left(2c_{2,i}{+}\frac{\z}{S_i}\left(\overline{\delta}_i {+}\underline{\delta}_i\right)\right) }\right]\!{,}~} \label{Eq: chi- def AB SW}\\
    \IEEEeqnarraymulticol{4}{l}{\chi^\text{+} \hspace{0.25\jot} = \frac{1}{|\mathcal{G}|} \!\sum_{i} \!\left[{{-} \nu \lambda_i {+} \left(\Sigma^\text{+}A_i\right)\left(2c_{2,i}{+}\frac{\z}{S_i}\left(\overline{\delta}_i {+}\underline{\delta}_i\right)\right) }\right]{,}} \label{Eq: chi+ def  AB SW}
    \end{IEEEeqnarray}
where $\Sigma^{-}$ and $\Sigma^{+}$ are the rows of covariance matrix $\Sigma$ related to the asymmetric forecast errors ($\Om^\text{-}$ and $\Om^\text{+}$).
Additionally, $S_i {=} \sqrt{A_i^{\!\top} \Sigma A_i} {=}\AZ$.
Note that the value of $S_i$ in the \OPFABSW{} model and in the derived balancing reserve policies presented below does not equal to the value of $S_i$ in the \OPFABNN{} model. Even though the expression for its calculation is the same, the way in which the vector $A_i$ and the matrix $\Sigma$ are formed depends on the balancing reserve policy itself.
Notably, the A-LBPs become dependent on the aggregated asymmetric uncertainty in the system and as such are equal for all the RES nodes, i.e., the generators are paid to balance the system imbalance independently of the sources contributing to this imbalance. Hence, relative to \OPFABNN{} model, a system-wide balancing approach does not capture a one-to-one relationship between individual sources of uncertainty and reserve providers, which causes inefficient price signaling. 

\subsection{Symmetric Node-to-Node Balancing -- \OPFSBNN{} \label{Sec: SB N2N}}
Another modification of \OPFABNN{} is derived under the assumption that controllable generators provide balancing reserve using nodal but symmetric participation factors ($\alpha_{iu}$) for symmetric forecast errors ($\om_u{=}\om_u^\text{-} + \om_u^\text{+}$) as in  \cite{Bienstock2018,Roald2016}. Under these assumptions, we obtain  $\alpha_{iu} {=} \alpha_{iu}^\text{-} {=} \alpha_{iu}^\text{+}$ and $A_i{=}\left[\alpha_{i1},\alpha_{i2},...,\alpha_{i|\mathcal{U}|} \right]^{\top}$. Random variable $\om$ is the vector of nodal forecast errors and $\Sigma$ is the covariance matrix of the nodal forecast errors.
\begin{IEEEeqnarray}{lrCll}
    \noalign{\textit{1) Generation Output:}
    Using these assumptions, \eqref{Eq: p AB N2N}, \eqref{Eq: Sufficiente a_iu AB}, \eqref{Eq: GCost AB N2N} can be replaced with:
    \vspace{-4\jot}} \IEEEnonumber\\
    &\IEEEeqnarraymulticol{3}{l}{\p_{i} = p_i - \!\textstyle{\sum_u}\left(\alpha_{iu}\om_u\right),} &\forall i \IEEEyesnumber\IEEEyessubnumber* \label{Eq: Total generation}\\
    &\IEEEeqnarraymulticol{3}{l}{{C_i}
    =   c_{i}{\left(p_{i}\right)} + c_{2,i}S_i,} &\forall i, \IEEEeqnarraynumspace \label{Eq: Expected SB NN}\\
    \left(\chi_u\right){:}\quad &
    \IEEEeqnarraymulticol{3}{l}{\hspace{0.5\jot} 1 \hspace{\jot}= \sum_u \alpha_{iu},} & \forall u.\label{eq:bAdN2n}\\
    \noalign{\noindent
    Expression \eqref{Eq: Expected SB NN} simplifies the expected balancing cost, which cancels out terms due to the symmetric forecast errors. 
    
    \vspace{2\jot}\textit{2) Probabilistic Capacity Limits:} Similarly to the expected cost in \eqref{Eq: Expected SB NN}, the capacity limits do not consider the asymmetric nodal RES forecast errors, which makes it possible to recast them as second-order conic constraints, which only depend on the variances of the nodal forecast errors:
    \vspace{-3\jot}} \IEEEnonumber\\
    \left(\overline{\delta}_{i}\right){:} \quad 
    & p_i + \z S_i &\leq& \hspace{2.5\jot} \overline{P}_i, &\forall i
    \label{Eq: Gen CC SB SOC NN}
    \label{Eq: Gen CC SB SOC NN Max}\\
    \left(\underline{\delta}_{i}\right){:} \quad 
    & \hspace{-2\jot} {-}p_i + \z S_i &\leq& {-}\underline{P}_i,&\forall i\label{Eq: Gen CC SB SOC NN Min}\\
    \noalign{\vspace{2\jot}\textit{3) Energy and Balancing Reserve Prices:} 
    Following  the  procedure  described  in  Proposition  \ref{Prop: AB N2N},  we derive the following prices: 
    \vspace{-3\jot}} \IEEEnonumber\\
    \IEEEeqnarraymulticol{4}{l}{\lambda_i = 2c_{2,i}p_i + c_{1,i} +\overline{\delta}_{i} - \underline{\delta}_{i},} & \forall i \label{Eq: KKT lambda_i SB NN}\\
    \IEEEeqnarraymulticol{4}{l}{
    \chi_u = \frac{1}{|\mathcal{G}|} \sum_{i} \left( \left(\Sigma_uA_i\right)\left( 2c_{2,i}+\frac{\z}{S_i}{\left(\overline{\delta}_{i} + \underline{\delta}_{i}\right)} \right)\right), ~} & \forall u, \IEEEeqnarraynumspace\label{Eq: KKT chi_u SB NN},
\end{IEEEeqnarray}
where $\Sigma_{u}$ is the row of the covariance matrix $\Sigma$ related to the nodal forecast error $\om_u$.
Thus, these LMPs do not change with the forecasted RES generation and the balancing reserve has only one price $\chi_u$, which is insensitive to the direction of the RES forecast errors and the LMPs. Note that the results in \eqref{Eq: KKT lambda_i SB NN} and \eqref{Eq: KKT chi_u SB NN} are similar those from \cite{Mieth2019,Dvorkin2020,Mieth2020}. 

\subsection{System-wide Symmetric Balancing -- \OPFSBSW{} \label{Sec: SB}}
The most studied chance-constrained electricity market design has   system-wide, symmetric RES forecast errors and symmetric generation response firstly introduced in \cite{Kuang2018} and then applied in \cite{Dvorkin2020,Mieth2019,Mieth2020}.
This case can be modeled within the proposed \OPFABNN{} model, if  $\Om$ is the system-wide RES forecast error, i.e., $\boldsymbol{\Omega} {=} \!\sum_u \boldsymbol{\omega}_u$,  $A_i^\text{-}{=}A_i^\text{+}{=}\left[\alpha_i, \forall u\right]$ and $S_i{=}S_i{=}\alpha_is$,  where $s^2$ is the total sum over the covariance matrix, $s^2 = e \Sigma e^T$.
\begin{IEEEeqnarray}{lrCll}
    \noalign{\textit{1) Generation Output:}
    Using these assumptions, \eqref{Eq: p AB N2N}, \eqref{Eq: Sufficiente a_iu AB}, \eqref{Eq: Gen c_i}, \eqref{Eq: GCost AB N2N} can be replaced with:
    \vspace{-4\jot}} \IEEEnonumber\\
    & \IEEEeqnarraymulticol{3}{l}{C_i = c^\text{}_i{\left(p_{i}\right)} + c_{2,i}^\text{}\alpha_i^2s^2,} &\forall i. \IEEEeqnarraynumspace\IEEEyesnumber\IEEEyessubnumber* \label{Eq: Exp Sb cost} \\
    &\IEEEeqnarraymulticol{3}{l}{\boldsymbol{p}_{i} = p_i - \alpha_i\boldsymbol{\Omega},} & \forall i \label{Eq: Total generation S}\\
    \left(\chi\right){:}\qquad&\IEEEeqnarraymulticol{3}{l}{\hspace{0.5\jot} 1 \hspace{\jot} = \sum_i\alpha_{i},} \label{Eq: Total Uncertainty Imbalance}\\
    \noalign{\vspace{2\jot}\textit{2) Probabilistic Capacity Limits:}
    Similarly to \eqref{Eq: Gen CC SB SOC NN Max}-\eqref{Eq: Gen CC SB SOC NN Min}, the capacity limits are reduced to the linear constraints:
    \vspace{-3\jot}} \IEEEnonumber\\
    \left(\overline{\delta}_{i}\right){:} \quad
    &p_i + \z \alpha_i s &\leq& \hspace{2.5\jot} \overline{P}_i, &\forall i\\
    \left(\underline{\delta}_{i}\right){:} ~~  \quad
    &\hspace{-2\jot} {-}p_i + \z \alpha_i s &\leq& {-}\underline{P}_i, \hspace{21\jot}&\forall i.\\
    \noalign{\vspace{2\jot}\textit{3) Energy and Balancing Reserve Prices:}
        Following  the  procedure  described  in  Proposition  2,  we derive the following prices: 
    \vspace{-3\jot}} \IEEEnonumber\\
    \IEEEeqnarraymulticol{4}{l}{\lambda_i = 2c_{2,i}p_i + c_{1,i} +\overline{\delta}_{i} - \underline{\delta}_{i},} & \forall i\label{Eq: KKT lambda_i SB}\\
    \IEEEeqnarraymulticol{4}{l}{
    \chi ~= \frac{1}{|\mathcal{G}|} \!\sum_{i} \!\left( {2c_{2,i}s^2\alpha_i} + 
        {\z s}
        {\left(\overline{\delta}_{i} + \underline{\delta}_{i}\right)} \right).} \IEEEeqnarraynumspace\label{Eq: KKT chi_u SB}
\end{IEEEeqnarray}

%

%
\begin{table*}[h]
\setlength\tabcolsep{3pt}
\caption{Variables for CC-OPF with balancing reserve \label{Tab: Gen values}}
   \centering
\resizebox{\linewidth}{!}{
\begin{tabular}{c|c|>{\centering\arraybackslash}m{1.1cm} >{\centering\arraybackslash}m{1.1cm}|c|cc}
\hline
\multicolumn{2}{c|}{\multirow{2}{*}{Balancing Type}}& \multicolumn{2}{c|}{Participation Factor} &Energy Price & \multicolumn{2}{c}{Balancing Price}\\
\cline{3-7}
\multicolumn{2}{c|}{} & $\alpha_{iu}^\text{-}$ & $\alpha_{iu}^\text{+}$ & $\lambda_i$ & $\chi_u^\text{-}$ & $\chi_u^\text{+}$ 
\\
\hline
&&&&&&\\[-1ex]
\multirow{2}{*}{Asymmetric} & Node-to-Node (OPF-N2N-AB) & $\alpha_{iu}^\text{-}$ & $\alpha_{iu}^\text{+}$ &  \multirow{2}{*}{$c_{2,i}\left(p_i-M \cdot A_i\right) {+} c_{1,i} {+} \overline{\delta}_{i} {-} \underline{\delta}_{i}$} & $\frac{1}{|\mathcal{G}|} \!\sum_{i} \!\left[{ \mu_u \lambda_i {+} \left(\SAm\right)\left(2c_{2,i}{+}\frac{\z}{S_i}\left(\overline{\delta}_i {+}\underline{\delta}_i\right)\right) }\right]$ & $\frac{1}{|\mathcal{G}|} \!\sum_{i} \!\left[{ {-}\mu_u \lambda_i {+} \left(\SAp\right)\left(2c_{2,i}{+}\frac{\z}{S_i}\left(\overline{\delta}_i {+}\underline{\delta}_i\right)\right) }\right]$
\\[2ex]
%
& System-wide (OPF-SW-AB)& $\alpha_{i}^\text{-}$ & $\alpha_{i}^\text{+}$ & & \multicolumn{2}{c}{$\frac{1}{|\mathcal{G}|} \sum_{i} \left( \left(\Sigma_uA_i\right)\left( 2c_{2,i}+\frac{\z}{S_i}{\left(\overline{\delta}_{i} + \underline{\delta}_{i}\right)} \right)\right)$}
\\[2ex]
\hline
&&&&&&\\[-1ex]
\multirow{2}{*}{Symmetric} & Node-to-Node (OPF-N2N-SB) & $\alpha_{iu}$ & $\alpha_{iu}$ & \multirow{2}{*}{$2c_{2,i}p_i + c_{1,i} +\overline{\delta}_{i} - \underline{\delta}_{i}$} & $\frac{1}{|\mathcal{G}|} \!\sum_{i} \!\left[{ \nu \lambda_i {+} \left(\Sigma^\text{-}A_i\right)\left(2c_{2,i}{+}\frac{\z}{S_i}\left(\overline{\delta}_i {+}\underline{\delta}_i\right)\right) }\right]$ & $\frac{1}{|\mathcal{G}|} \!\sum_{i} \!\left[{ {-}\nu \lambda_i {+} \left(\Sigma^\text{+}A_i\right)\left(2c_{2,i}{+}\frac{\z}{S_i}\left(\overline{\delta}_i {+}\underline{\delta}_i\right)\right) }\right]$
\\[2ex]
%
& System-wide (OPF-SW-SB)& $\alpha_{i}$ & $\alpha_{i}$ & & \multicolumn{2}{c}{$\frac{1}{|\mathcal{G}|} \!\sum_{i} \!\left( {2c_{2,i}s^2\alpha_i} + 
        {\z s}
        {\left(\overline{\delta}_{i} + \underline{\delta}_{i}\right)} \right)$}\\[2ex]
\hline
\hline
\end{tabular}
}
\end{table*}

\section{Discussion and Insights on  Locational Energy and Balancing Reserve Prices}
\subsection{Symmetric vs  Asymmetric Energy Pricing}
In the asymmetric balancing policies,  energy price $\lambda_i$ reflects not only incremental production costs through coefficients $c_{2,i}$ and $c_{1,i}$, and the generation capacity limits through duals $\overline{\delta}_i$ and $\underline{\delta}_i$, but also the expected generation level ($p_i$) and the expected real-time balancing power ($M{\cdot}A_i$), i.e., $p_i - M{\cdot}A_i$.
Therefore, the expected asymmetric balancing provision affects the LMPs by reducing or increasing the expected available generation at the generation nodes.
The influence of the expected balancing reserve provision on the LMPs does not occur in symmetric balancing markets, because the expected balancing reserve provision equals zero.
\subsection{Symmetric vs  Asymmetric Balancing Reserve Pricing}
Unlike symmetric LBPs, the A-LBPs include a term that relates  LMPs $\lambda_i$ to the nodal expected asymmetric forecast error $\mu_u$ in \eqref{Eq: KKT chi_u AB NN m}--\eqref{Eq: KKT chi_u AB NN p} and to the expected system-wide forecast error $\nu$ in \eqref{Eq: chi- def AB SW}--\eqref{Eq: chi+ def AB SW}.
Therefore, the A-LBPs reflect the expected balancing need at node $u$ and the LMP of the balancing generator. 

\subsection{Nodal vs System-Wide Balancing Reserve Pricing}
\label{sec:relSymPrices}
To analyze the effect of node-to-node balancing reserve policies, we evaluate the relationship between the LBPs and the system-wide balancing reserve prices.
To simplify such comparison, we will focus on the symmetric balancing reserve policies, because the asymmetric balancing analysis leads to more complicated expressions.

Under symmetric balancing reserve policies,  
the system-wide and node-to-node participation factors are equal, i.e.,  ${\alpha_{i} {=} \alpha_{iu}}, \text{ and } \alpha_i\Om {=} \!\sum_u \!\left(\alpha_{iu}\om_u\right);~ \forall i,u$. Then, it follows that:
\begin{IEEEeqnarray}{Cl}
    S_i = \soc{A_i\Sigma^{\nicefrac{1}{2}}} = \alpha_i s , \qquad &\forall i. \label{Eq: alphas equiv SB}
\end{IEEEeqnarray}
It is possible to derive a relationship $\beta_u$ between the LBPs and the system-wide balancing reserve price, where:
\begin{IEEEeqnarray}{rCll}
    \beta_u &=& \frac{\chi_{u}}{\chi}
    , \qquad &\forall u. \IEEEeqnarraynumspace \label{Eq: Beta u def}
\end{IEEEeqnarray}
An analytical expression for \eqref{Eq: Beta u def} can be derived by factorizing \eqref{Eq: KKT chi_u SB NN} in the following way:
\begin{subequations}
\begin{IEEEeqnarray}{rCll}
    \chi_{u} &{=}& \frac{1}{|\mathcal{G}|} \sum_{i} \left( \left(\Sigma_uA_i\right)\left( 2c_{2,i}+\frac{\z}{S_i}{\left(\overline{\delta}_{i} + \underline{\delta}_{i}\right)} \right)\right)
    , \quad    &{\forall u}\\
    \noalign{\noindent using identity \eqref{Eq: alphas equiv SB} and since $\alpha_{iv} = \alpha_i, \forall i,v$ \vspace{2\jot}}
    \chi_{u} &{=}& \frac{1}{|\mathcal{G}|}\!\sum_i\!\left[
    \!\sum_{v} \!\left( \alpha_{i} s_{uv}\right) \!\left(
    {2c_{2,i}  {+} {\z}
    \frac{\overline{\delta}_{i} {+} \underline{\delta}_{i}}{\alpha_i s}
    }\right)\right],
     &\forall u \IEEEeqnarraynumspace \\
    &{=}& \frac{1}{|\mathcal{G}|}\!\sum_i\!\left[\frac{\!\sum_{v} s_{uv}}{s^2} \!\left(
    {2c_{2,i}{s^2\alpha_i}  {+} {\z s}
    \left({\overline{\delta}_{i} {+} \underline{\delta}_{i}}\right)
    }\right)\!\right],
    &\forall u \IEEEeqnarraynumspace \\
   &{=}& \frac{\!\sum_{v} s_{uv}}{ s^2} \!\left[\frac{1}{|\mathcal{G}|}\!\sum_i \!\left(
    {2c_{2,i}{s^2\alpha_i}  {+} {\z s}
    \left({\overline{\delta}_{i} {+} \underline{\delta}_{i}}\right)
    }\right)\right]{,}
    &\forall u \IEEEeqnarraynumspace \\
    \noalign{\noindent where the term inside the square brackets is equal to  balancing reserve price in $\chi$ \eqref{Eq: KKT chi_u SB}. 
    Using $s^2=e\Sigma e = \sum_{v,w}\sigma_{vw}$ yields:
    \vspace{2\jot}}
    \IEEEeqnarraymulticol{3}{c}{\chi_{u} {=} \frac{\!\sum_{v}\sigma_{uv}}{\!\sum_{v,w}\sigma_{vw}} \chi
    {,}} &\forall u. \IEEEeqnarraynumspace\label{eq:chiU}\\
    \noalign{\noindent Finally,  rearranging the terms in \eqref{eq:chiU} to match     \eqref{Eq: Beta u def} leads to:
    \vspace{\jot}}
    \IEEEeqnarraymulticol{3}{c}{\beta_{u} = \frac{\!\sum_{v}\sigma_{uv}}{\!\sum_{v,w}\sigma_{vw}}
    {,}} &\forall u. \IEEEeqnarraynumspace \label{Eq: beta u analytical}
\end{IEEEeqnarray}
\end{subequations}
Parameter $\beta_u$ does not depend on system parameters, primal and dual variables. Thus, the nodal balancing reserve price depends on the effect that the RES at node $u$ has on the system-wide variance, i.e., the greater the introduced uncertainty by RES at node $u$, the greater the price that it must pay for deviations from its forecast. Further, it follows that  $\sum_u \beta_u = 1$, which makes it possible to interpret $\beta_u$ as the contribution of the RES at node $u$ to the overall system balancing reserve price.

Since the LBPs $\chi_u$ carry information on the marginal balancing cost that the RES at node $u$ adds to the system, they can be used to incentivize system-beneficial RES deployment strategies. 
For example, a RES investor would prefer to install their RES project at a node with low, or ideally negative, correlation with the existing RES, so payments for introducing uncertainty and increasing balancing reserve requirements will be lower, possibly leading to a neutral or even negative value of $\beta_u$. The balancing reserve prices derived from the proposed framework highlight the value that the system as a whole derives from complementary RES generation, i.e., RES generators with a negative correlation to other RES,  since a RES that compensates the imbalance of others provides load balancing reserve support and as such its balancing costs are lower.
\section{Analysis of RES Costs}
Typically, RES resources submit zero-price bids to electricity market operators, reflecting their zero or near-zero marginal production cost.
However, they inflict a non-zero marginal cost on the system by driving the need for procuring additional balancing reserves.
This section compares two approaches to quantify these costs by (i) using the balancing reserve payments to controllable generation resources and by (ii) using the marginal balancing cost. 
We show that both approaches can be equivalent.
\subsection{Compensating Balancing Cost}
We denote the the balancing cost incurred by the system for compensating a nodal power imbalanced caused by the RES at node $u$ as $C_u^\alpha$ and define it as the sum of payments to reserve providers:
\begin{IEEEeqnarray}{rClll}
    \IEEEeqnarraymulticol{4}{C}{C_u^\alpha =  \chi_u \!\sum_i \alpha_{iu} , \qquad} & \forall u\label{Eq: C_u def} \IEEEyesnumber\IEEEyessubnumber*\\
    \noalign{\noindent and since $\!\sum_i \alpha_{iu}{=}1$ we obtain: 
    \vspace{2\jot}}
    C_u^\alpha &{=}& \chi_u {=}
    {\frac{1}{|\mathcal{G}|}} \!\sum_{i}&
    \left[ \frac{\sum_{v} \alpha_{iv} \sigma_{uv}}{S_i}
    \!\left[2c_{2,i} S_i
    {+} {\z}
    {\!\left(\overline{\delta}_{i} {+} \underline{\delta}_{i}\right)}
    \right]\right]{,} & \forall u. \IEEEeqnarraynumspace \label{Eq: C_u def ext}
\end{IEEEeqnarray}

The balancing reserve compensation cost $C_u^\alpha$ paid by the RES at node $u$ depends on the uncertainty it introduces to the system ($\sigma_{u,u}$) and its correlation with other nodes $v$ ($\sigma_{u,v}$). If the RES at $u$ does not introduce uncertainty to the system, then $C_u^\alpha{=}0$.

\subsection{RES Costs as Marginal Balancing Costs}
Since the RES are modelled with a zero marginal operating cost, their revenue can be calculated as:
\begin{IEEEeqnarray}{C}
    R_u = -\frac{\partial \mathcal{L_\text{SN}}}{\partial \text{w}_u} \cdot \text{w}_u, \qquad \forall u \label{Eq: Profit u def}
\end{IEEEeqnarray}
where \(\mathcal{L_\text{SN}} \) is the Lagrangian function for the \OPFSBNN{} model under the node-to-node symmetric balancing reserve policy presented in Section \ref{Sec: SB N2N}. 
If the standard deviation $\sigma_u$ of $\om_u$ is proportional to its forecast injection, i.e. $\sigma_u {=} \kappa_u \text{w}_u$, and  $\zeta_{uv}$ is the correlation coefficient between uncertain injections at nodes $u$ and $v$ \cite{Fang2019}.
Then, the marginal revenue for the RES resource can be computed as:
\begin{IEEEeqnarray}{lr}
    \IEEEeqnarraymulticol{2}{l}{-\frac{\partial \mathcal{L_\text{SN}}}{\partial \text{w}_u} {=} \lambda_u {-}{\!\sum_i \!\left(\frac{\!\sum_v \alpha_{iu}\alpha_{iv}\zeta_{uv}\kappa_u\sigma_v}{S_i} 
    \!\left[ 2c_{2,i}S_i {+} \z\!\left(\overline{\delta}_{i} {+} \underline{\delta}_{i} \right) \right]\!\right)\!,}}\IEEEnonumber\\
    &\forall u \label{Eq: Profit u} \IEEEeqnarraynumspace
\end{IEEEeqnarray}
where the positive part of expression \eqref{Eq: Profit u} can be interpreted as the active power price paid to the RES at node $u$, $\pi_u^\text{p}$, and the negative one as the cost associated with its stochastic generation, $c_u$. 
If the the RES at node $u$ introduces no uncertainty to the system, i.e., ${\kappa_u {=} 0}$, then $c_u{=}0$.
The operation cost of the the RES at node $u$ is then given by:
\begin{IEEEeqnarray}{c}
    C^\text{w}_u = c_u \text{w}_u, \qquad \forall u. \label{Cp_u}
\end{IEEEeqnarray}
\subsection{Cost Equivalence between the Compensating and Marginal Balancing Costs}
The equivalence between $C_u^\alpha$ and $C_u^\text{w}$ can be  demonstrated by examining \eqref{Cp_u}:
\begin{IEEEeqnarray}{rClll}
    C^\text{w}_u &{=}& c_u \text{w}_u \IEEEnonumber\\
    &{=}& \IEEEeqnarraymulticol{2}{l}{\text{w}_u {\!\sum_i \!\left(\frac{\!\sum_v \alpha_{iu}\alpha_{iv}\zeta_{uv}\kappa_u\sigma_v}{S_i} 
    \!\left[ 2c_{2,i}S_i {+} \z\!\left(\overline{\delta}_{i} {+} \underline{\delta}_{i} \right) \right]\!\right)\!,}} \IEEEnonumber\\
    &&& \forall u \IEEEyesnumber\IEEEyessubnumber*\\
    &{=}& \IEEEeqnarraymulticol{2}{l}{{\!\sum_i \!\left(\frac{\!\sum_v \alpha_{iu}\alpha_{iv}\zeta_{uv}\kappa_u\text{w}_u\sigma_v}{S_i}
    \!\left[ 2c_{2,i}S_i {+} \z\!\left(\overline{\delta}_{i} {+} \underline{\delta}_{i} \right) \right]\!\right)\!,}} \IEEEnonumber\\
    &&& \forall u. \\
    \noalign{\noindent Next, using $\sigma_u {=} \kappa_u \text{w}_u$ and $\sigma_{uv} {=} \zeta_{uv}\sigma_u\sigma_v$, it follows:
    \vspace{2\jot}}
    C^\text{w}_u
    &{=}& \IEEEeqnarraymulticol{1}{l}{\!\sum_i  \left(\frac{ \sum_v \alpha_{iu}\alpha_{iv}\sigma_{uv}}{S_i}
    \left[ 2c_{2,i}S_i {+} \z \left(\overline{\delta}_{i} {+} \underline{\delta}_{i} \right) \right] \right) ,} 
    & \forall u \IEEEeqnarraynumspace\\
    &{=}& \IEEEeqnarraymulticol{1}{l}{\alpha_{iu}\!\sum_i  \left(\frac{ \sum_v \alpha_{iv}\sigma_{uv}}{S_i}
    \left[ 2c_{2,i}S_i {+} \z \left(\overline{\delta}_{i} {+} \underline{\delta}_{i} \right) \right] \right) ,}& \forall u \IEEEeqnarraynumspace\\
    \noalign{\noindent where replacing $\chi_{u}$ results in:
    \vspace{2\jot}}
    \IEEEeqnarraymulticol{3}{C}{C^\text{w}_u = \chi_{u} \!\sum_i \alpha_{iu} = C^\alpha_{u},} &\forall u. \IEEEeqnarraynumspace \label{Eq: equivalence C_u}
\end{IEEEeqnarray}
As shown in  Eq.~\eqref{Eq: equivalence C_u}, the total operating cost of the RES resources is equivalent when computed from balancing reserve payments $C^\alpha_{u}$ and  marginal costs $C^\text{w}_{u}$. However, the latter holds only if the total operating cost of the RES resources is computed under the assumptions that the variances, \(\sigma_{u,v}\), are  proportional to the injected power and fixed  node-to-node correlation coefficients, i.e., $\sigma_u {=} \kappa_u\text{w}_u$.  If these assumptions do not hold, the total operating cost of the RES resources can only be computed via the balancing reserve payments in \eqref{Eq: C_u def ext}.

\section{Numerical Experiments \label{Sec: Case Study}}
In this section, the proposed \OPFABNN~ model and its modifications from Section~\ref{Sec:modifications} are evaluated on a modified IEEE 118-bus case \cite{Mieth2020}. The experiments were conducted using Julia 1.53, JuMP 0.21.6 and Mosek \cite{Mosek} on an Intel i7-1165G7 processor clocked at 2.80GHz with 16 GB of RAM.
The code and data supplement can be downloaded from \cite{AntonRepo}.

\subsection{Test Case and Stochastic Data}

For this case study, we add 11 utility-scale wind farms to the IEEE-118 bus test system as the RES resources. Table~\ref{tab:parameters} indicates their placement in the system. 
Forecasts and forecast error distributions are emulated from real data of onshore wind farms in different European regions collected between January $\text{1}^{\text{st}}$ and October $\text{1}^{\text{st}}$ of 2020 from  \cite{ENTSO-E}. 
Table~\ref{tab:parameters} summarizes the normalized wind power data. 

To study the RES impact at different penetration levels, we create four renewable penetration scenarios by scaling the nominal installed capacity of each wind farm by 50\%, 100\%, 200\%, and 400\%. 
In these scenarios, the available generation from all wind farms covers 4.8\%, 9.7\%, 19.4\%, and 38.8\% of total demand, respectively.

\begin{table}[b]
\caption{Normalized Forecasts $w_{u}$
mean and standard deviations at nodes $u$}
\label{tab:parameters}
\begin{center}
\vspace{-1em}
\resizebox{\columnwidth}{!}{%
{\def\arraystretch{1.2}\tabcolsep=4pt
\begin{tabular}{|c|ccccccccccc|}
\hline
$u$ & 3 & 8 & 11 & 20 & 24 & 26 & 31 & 38 & 43 & 49 & 53\\
\hline
$w_{u}$ & 0.112 & 0.801 & 0.61 & 0.086 & 0.142 & 0.0 & 0.056 & 0.137 & 0.353 & 0.207 & 0.305\\
\hline
$\sigma_{u}$ & 0.1 & 0.26 & 0.22 & 0.25 & 0.22 & 0.19 & 0.14 & 0.44 & 0.17 & 0.13 & 0.07\\
\hline
$\mu_{u}$ & 0.1544 & 0.4101 & 0.3539 & 0.3962 & 0.3483 & 0.3029 & 0.2267 & 0.7016 & 0.2791 & 0.2113 & 0.1106\\
$\sigma^{\text{-}}_{u}$ & 0.06 & 0.155 & 0.142 & 0.116 & 0.121 & 0.097 & 0.087 & 0.279 & 0.103 & 0.079 & 0.039\\
$\sigma^{\text{+}}_{u}$ & 0.056 & 0.155 & 0.126 & 0.19 & 0.142 & 0.134 & 0.084 & 0.251 & 0.108 & 0.081 & 0.045\\
\hline
\end{tabular}}}
\end{center}
\end{table}

Every wind farm was assigned an independent normal distribution  $\om_{u}\sim\mathcal{N}(0, \sigma_{u})$ with $\sigma_{u}$ being estimated from the data in \cite{ENTSO-E}. 
For asymmetric balancing reserve policies, we generate  truncated normal distributions $\om_{u}^{\text{+}}$ and $\om_{u}^{\text{+}}$ from  $\om_{u}$ with the following parameters:
\begin{IEEEeqnarray}{c}
    \sigma_{u}^{\pm} = \sigma_{u}\sqrt{\frac{2\pi - 4}{2\pi}},\quad \mu_{u}^{\pm} = \sigma_{u}\sqrt{\frac{2}{\pi}} \label{eq:pmdist}
\end{IEEEeqnarray}

We set the confidence level of chance constraints at $(1{-}\epsilon) {=} 0.99$, obtaining $z_i = \sqrt{(1-\epsilon)/\epsilon}\approx 9.95$ with the Chebyshev inequality as in \cite{Dvorkin2020}.

\subsection{Solution Approach}
\label{sec:optimization}
Chance-constrained market clearing procedures as proposed in \cite{Kuang2018,Dvorkin2019,Mieth2019,Mieth2020} can be solved in a single or two-step
procedure using off-the-shelf solvers.   
However, solving \OPFABNN~and its  modifications shown in Section~\ref{Sec:modifications} requires dealing with bilinear terms in the objective function. 
Therefore, we use McCormick envelopes to convexify the bilinear terms in \eqref{Eq: Obj OPFABNN}. The obtained envelopes are then sequentially tightened around the  operating point produced at the previous iteration to improve accuracy.  
The proposed sequential solution approach based on using off-the-shelf solvers is described below:

\begin{itemize}
    \item \textit{Step 1.} Obtain a solution to a convex approximation of the original \OPFABNN~problem.
    \item \textit{Step 2.} Tighten the convex approximation around \textit{Step 1}'s solution by using a decreasing scalar factor to reduce the distance between its obtained operation point and the current bounds.
    
    \item \textit{Step 3.} Stop if the convex approximation is good enough. If not, update and return to \textit{Step 1} and repeat.
\end{itemize}

This solution approach has two advantages.
First, the convexified problem always yields a feasible solution to the \OPFABNN{} optimization, because the convexification only affects the objective function.
Second, an optimal solution $\hat{x}^{*}$ of the convexified problem provides an upper bound to the objective of the original non-convex problem, i.e., $C(x^*) \leq {C}(\hat{x}^{*})$, where $C(\cdot)$ is the objective function of \OPFABNN{}, and $x^*$ is an optimal solution to the original non-convex  \OPFABNN{} version.
On the other hand, the optimal objective value of the convexified problem $\hat{C}(\hat{x}^{*})$ provides a lower bound to the objective of the original non-convex problem, i.e., $\hat{C}(\hat{x}^{*}) \leq C(x^*)$.
Therefore, we can observe that $\hat{C}(\hat{x}^{*}) \leq C(x^*) \leq C(\hat{x}^{*})$. 
However, while $C(x^*)$ cannot be computed directly, both values $\hat{C}(\hat{x}^{*})$ and $C(\hat{x}^{*})$ can be obtained directly. 
Thus, the stopping criteria threshold from \textit{Step 3} is set by the approximation accuracy $100(\nicefrac{C(\hat{x}^{*})}{\hat{C}(\hat{x}^{*})}-1)$, representing an optimality gap. Once the approximation accuracy is below this threshold, we can derive the prices using the convexified \OPFABNN{} version.

\subsection{Comparison of Balancing Reserve Policies}
\label{sec:scenarios}

\begin{figure}[b]
\begin{center}
{\includegraphics[width = \linewidth]{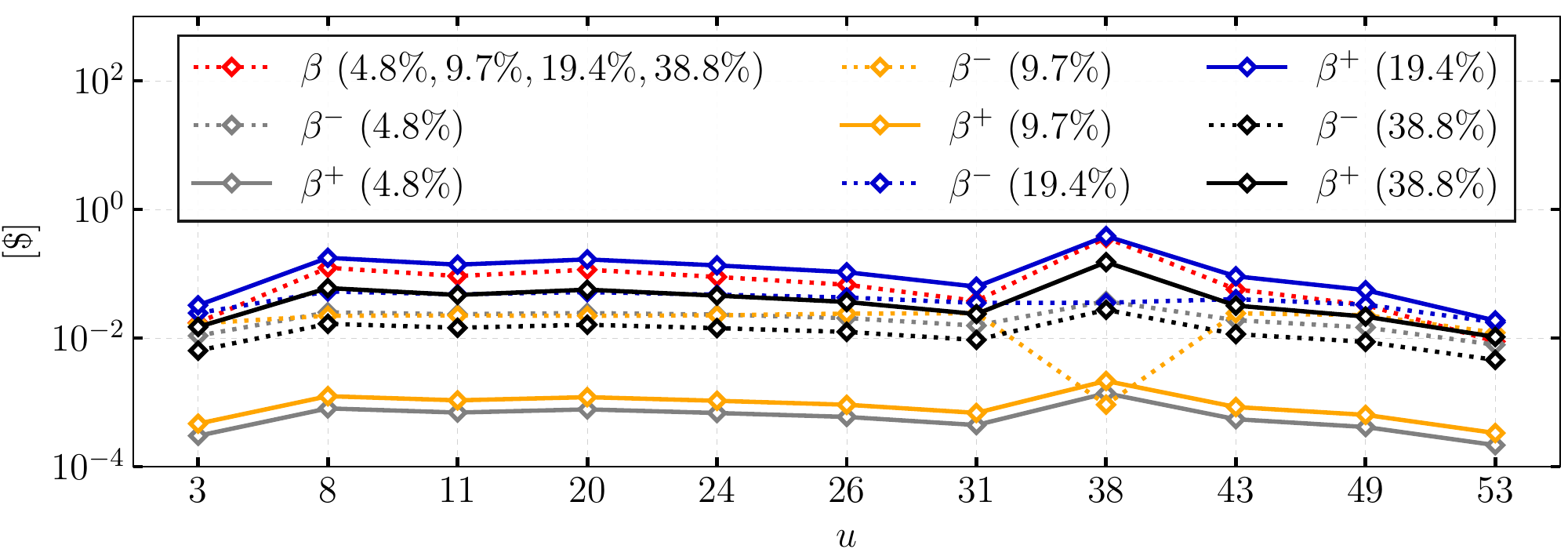}}\\
{\includegraphics[width = \linewidth]{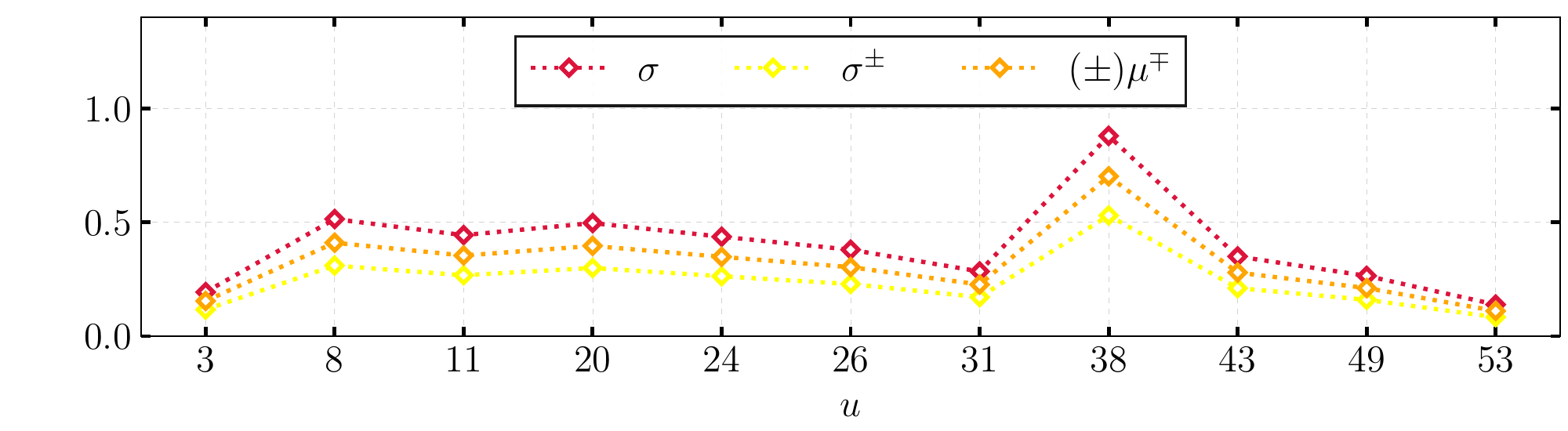}}
  \caption{Relationship between $\beta^{(\pm)}_u$ and stochastic generation parameters.\label{fig:betas}}
\end{center}
\end{figure}

\subsubsection{System-wide and node-to-node balancing}
\label{sec:sym}
Table \ref{tab:resultszeromean} summarizes the main results for our base scenario with a renewable penetration level of 9.7\%. 

Our results show that symmetric, zero-mean balancing yields the same objective value in both system-wide and node-to-node policies. 
Moreover, by inspecting dual values $\chi$ and $\chi_u$ we observe that $\sum_{u}\chi_{u} = \chi$ as shown analytically in Eq. \eqref{Eq: Beta u def} above. 
Furthermore, the dependence of $\chi_u$ on the variance $\sigma^2_u$ is illustrated in Figure \ref{fig:betas}. 
As such, a node-to-node balancing reserve policy provides better pricing signals by reflecting the influence the forecast error variance introduced by the RES has on the price paid for its balancing.

{\def\arraystretch{1.25}\tabcolsep=2pt
\begin{table}[b]
\caption{Computational Results for Markets with Symmetric Balancing Provision (SW, N2N)}
\label{tab:resultszeromean}
\begin{center}
\begin{tabular}{|c|c|c|}
\hline
& SW & N2N\\
\hline
$z^*~[\$/h]$ & 127020.2 & 127020.2\\
\hline
$\chi$ & 22.9 & -\\
\hline
\end{tabular}
\end{center}
\vspace{1ex}
    \resizebox{\columnwidth}{!}{%
    \begin{tabular}{|c|ccccccccccc|c|}
        \hline
        $u$ & 3 & 8 & 11 & 20 & 24 & 26 & 31 & 38 & 43 & 49 & 53 & $\sum$\\
        \hline
        $\chi_u$ & 0.4 & 2.822 & 2.102 & 2.634 & 2.036 & 1.539 & 0.862 & 8.258 & 1.307 & 0.749 & 0.205 & 22.9\\
        \hline
    \end{tabular}}
\end{table}}

\subsubsection{Asymmetric balancing provision with varying amounts of renewable penetration}
\begin{figure}[t]
\vspace{-1em}
\begin{center}
\subfloat[4.8\% renewable penetration]{\includegraphics[width = .5\linewidth]{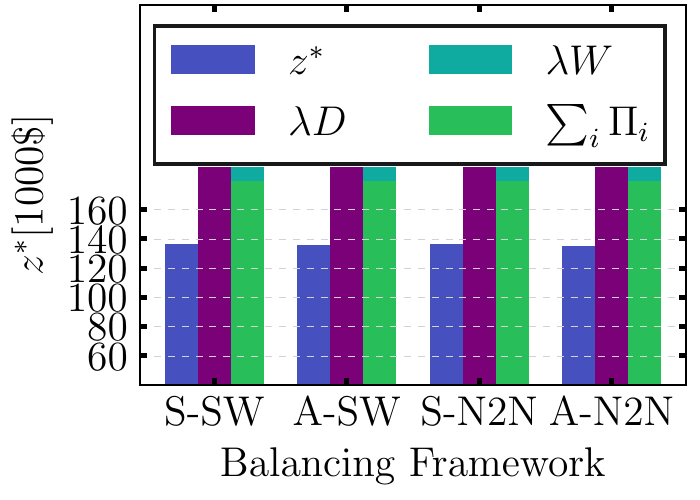}\label{fig:obj_5}}
\subfloat[9.7\% renewable penetration]{\includegraphics[width = .5\linewidth]{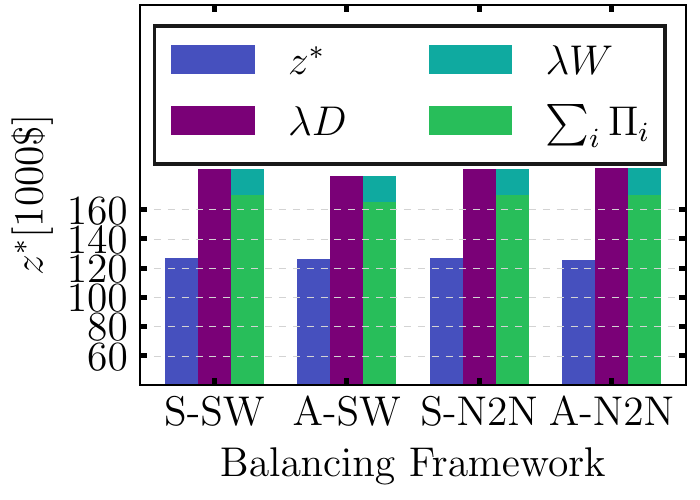}\label{fig:obj_10}}\\
\subfloat[19.4\% renewable penetration]{\includegraphics[width = .5\linewidth]{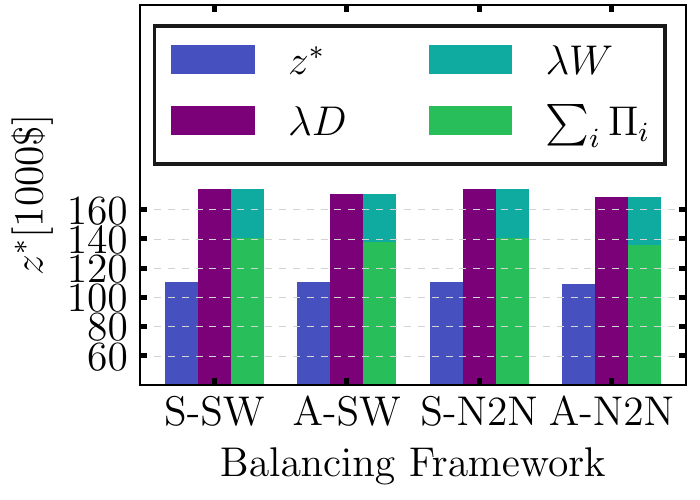}\label{fig:obj_20}}
\subfloat[38.8\% renewable penetration]{\includegraphics[width = .5\linewidth]{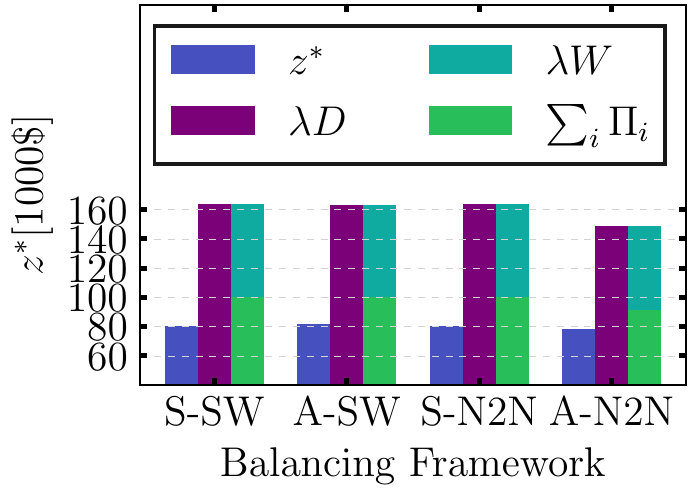}\label{fig:obj_40}}
  \caption{Optimization results under different balancing policies and renewable penetration. Where $z^*$ represents the objective value, $\lambda D$ consumer payments, $\lambda W$ wind payments and $\sum_{i}\Pi_{i}$ the generators' revenue.\label{fig:objectives}}
\end{center}
\end{figure}

Next, we study the asymmetric balancing reserve provision as described in Section~\ref{sec:optimization} with different renewable penetration levels.
To facilitate the comparison, the results in Figure \ref{fig:objectives} are presented with the results obtained using the symmetric balancing policy (as in Section \ref{sec:sym} above). 
Figure \ref{fig:objectives} contains the best objective value found in the sequential process described in Section~\ref{sec:optimization}, such that $C(\hat{x}^*) \approx C(x^*)$. Further, it shows the payments collected from consumers ($\lambda D$), payments made to wind farms ($\lambda W$) and conventional generators ($\sum_i\Pi_i$ ), respectively, for the different renewable penetration levels defined above. 

With an increasing penetration of wind generation, system cost given by objective value $C(x^*)$ are decreasing due to the greater availability of low-cost energy from RES. This leads to notable differences in payments made to conventional generators. However, for low renewable penetration rates, the total cost barely change between balancing reserve policies.

Our experiments confirm the revenue adequacy property of the prices, i.e., $\lambda D = \lambda W + \sum_i\Pi_i$.
Table \ref{tab:payments} shows consumer payments, generator revenues and generation costs using different balancing reserve policies.
The difference in consumer payments between asymmetric and symmetric balancing markets is defined as $\Delta\lambda D$.
The approximation error introduced by the convexification is captured in the term $err$.
As summarized in Figure \ref{fig:objectives} and Table \ref{tab:payments}, markets with asymmetric node-to-node balancing policies show that significant reductions in consumer payments are possible, representing up to 9\% of cost savings.

These results point to two observations. First, the solution procedure in  Section~\ref{sec:optimization} performs well for system-wide asymmetric balancing policies. Given a solution, the objective differs by less than 0.0026\% from the exact objective. Second, the proposed  node-to-node balancing reserve policies yield consistent reductions in consumer payments relative to the symmetric balancing reserve policy.

{\def\arraystretch{1.2}\tabcolsep=4pt
\begin{table}[t]
\caption{Consumer Payments $\lambda D$  and Approximation Errors by Scenario}
\label{tab:payments}
\begin{center}
\begin{tabular}{|c|c|cccc|}
\hline
\multicolumn{2}{|c|}{Scenario} & 4.8\% & 9.7\% & 29.4\% & 38.8\%\\
\hline
S-SW  & $\lambda D~[\$\backslash h]$ & 188965.12 & 188082.21 & 174123.79 & 164091.47\\
\hline
\multirow{ 3}{*}{A-SW} & $\lambda D~[\$\backslash h]$ & 188940.56 & 182724.69 & 170674.6 & 163494.45\\
& $\Delta\lambda D~[\%]$ & -0.013 & -2.85 & -1.98 & -0.364 \\ 
& $\mathrm{err}~[\%]$ & 4.0e-5 & 1.0e-4 & 3.8e-4 & 2.63e-3 \\ \hline
S-N2N & $\lambda D~[\$\backslash h]$ & 188965.16 & 188082.55 & 174123.41 & 164091.43\\
\hline
\multirow{ 3}{*}{A-N2N} & $\lambda D~[\$\backslash h]$ & 188869.01 & 188465.01 & 168506.01 & 148982.8\\
& $\Delta\lambda D~[\%]$ & -0.051& +0.2 & -3.23 & -9.21 \\ 
& $\mathrm{err}~[\%]$ & 0.70217 & 1.01088 & 1.15133 & 1.79429 \\
\hline
\end{tabular}
\end{center}
\end{table}}

\begin{figure}[t]
\begin{center}
{\includegraphics[width = \linewidth]{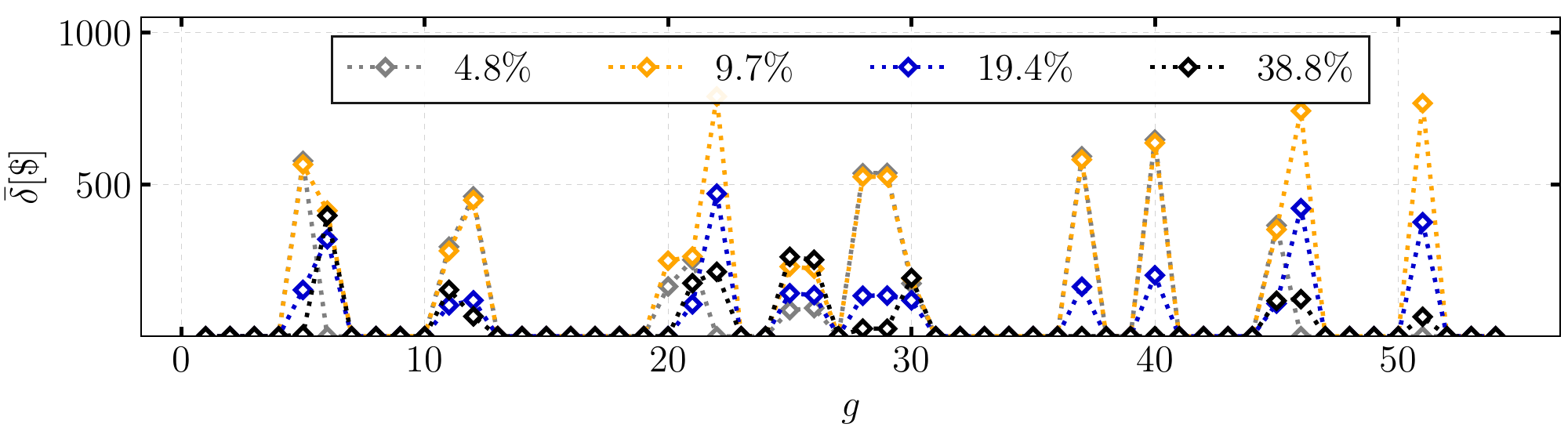}}\\
{\includegraphics[width = \linewidth]{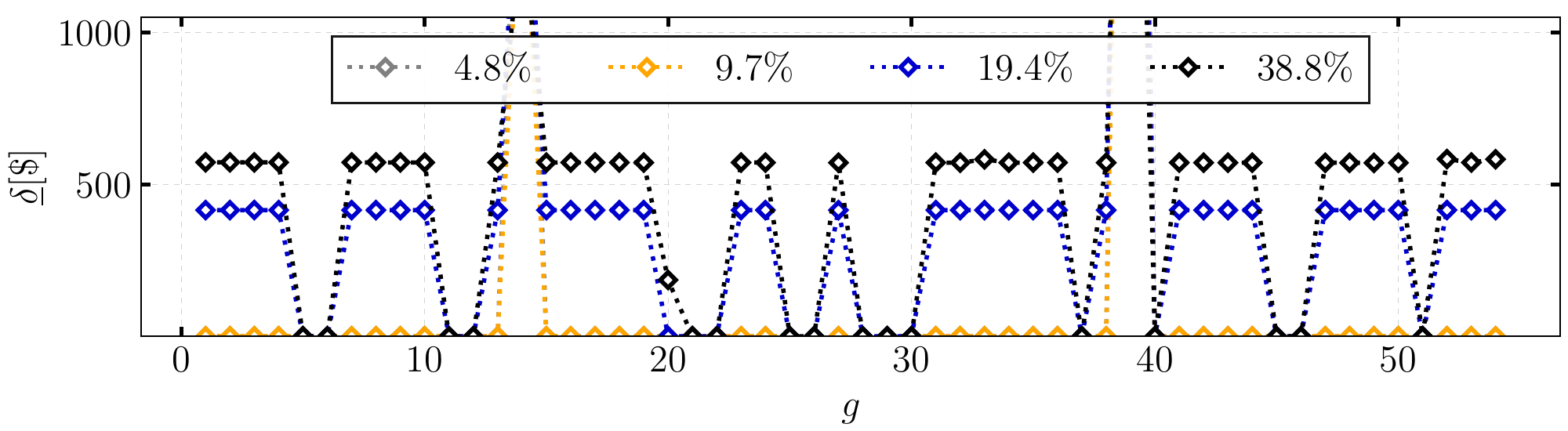}.}
  \caption{Dual values $\underline{\delta}$ and $\overline{\delta}$ for different renewable capacity levels.\label{fig:deltas}}
\end{center}
\end{figure}

Figure \ref{fig:deltas} provides insights as to how prices are formed under the proposed asymmetric node-to-node balancing policy. We can observe that the dual variables $\overline{\delta}$ are proportional to the renewable penetration level. This is due to fewer dispatchable generation capacity required, i.e., less expensive generators must be dispatched to clear the market. Simultaneously, the variance introduced into the system increases with greater renewable penetration level, which leads to tighter chance constraints \eqref{Eq: Pmin u AB} and \eqref{Eq: Pmax u AB}. Accordingly, more dispatchable generation capacity must be called for balancing reserve provision, which leads to an increase in dual variables and, subsequently, balancing reserve prices.

\section{Conclusion}
In this paper, we have described and analyzed the economic impact of four policies for the provision of balancing regulation for the chance-constrained electricity market: system-wide symmetric, node-to-node symmetric, system-wide asymmetric, and node-to-node asymmetric. 
Both the locational marginal electricity prices and the balancing regulation prices are derived and are proven to establish a competitive equilibrium under all the balancing policies.
The asymmetric node-to-node regulation provides the least balancing prices, while introducing a dependency between the  electricity prices and the generation participation in the balancing provision.
We have shown that asymmetric reserve provision facilitates more efficient energy and reserve pricing by co-aligning them with uncertainty information and by minimizing the operating cost, which in turn leads to lower demand charges, an important property to accommodate large  RES portfolios. 
Additionally, other balancing policies are proven to be particular cases of the asymmetric node-to-node regulation one, thus yielding greater operational costs.


\bibliography{Reserves_Pricing}
\bibliographystyle{IEEEtran}


\appendix
\numberwithin{equation}{subsection}
\subsection{Expected Cost Derivation: Asymmetric Node-to-Node Balancing}\label{App: Asymmetric cost N2N} 
Generator $i$'s cost can be expressed in terms of nodal directional imbalances by
\begin{IEEEeqnarray}{rCll}
    {c_i}{(\boldsymbol{p}_{i})} & = & {c_{2,i}} {{\left( \boldsymbol{p}_{i} \right)}^{2}} {+} c_{1,i} {\left( \p_{i} \right)} {+} c_{0,i}, 
    & \forall i \IEEEeqnarraynumspace \\
    \noalign{\noindent where
    \vspace{2\jot}}
    \boldsymbol{p}_{i} &=& p_{i} {-} \!\sum_u \left( \alpha^{-}_{iu} \om_u^- {+} \alpha^{+}_{iu}\om_u^+ \right), \quad &\forall i. \IEEEyessubnumber
\end{IEEEeqnarray}
\noindent Let the vectors $A_i{=}\left[\alpha_{i1}^-,\alpha_{i2}^-,...,\alpha_{i|U|}^-, \alpha_{i1}^+, \alpha_{i2}^+,..., \alpha_{i|U|}^+ \right]^{\top}$ , $\boldsymbol{\Omega}{=}\left[\om_{1}^-,\om_{2}^-,...,\om_{|U|}^-,\om_{1}^+,\om_{2}^+,...,\om_{|U|}^+ \right]$, then the expected generation costs can be calculated as:
\begin{IEEEeqnarray}{rCll}
    \Expected{c_i{(\p_i)}} & = & \Expected{c_i{(p_i{-}\Om \cdot A_i)}}
    \IEEEnonumber\\
    & = &
    c_i\left( p_{i} \right) + \Expected{\Y_{1,i}} + \Expected{\Y_{2,i}}, & \qquad \forall i \IEEEeqnarraynumspace \label{Eq: Exp AB Yp Ym N2N}
\end{IEEEeqnarray}
\noindent where
\begin{IEEEeqnarray}{rCll}
    \Y_{1,i} & = & -\left(2 c_{2,i} p_i {+} c_{1,i}\right) \X_i,
    & \quad \forall i \IEEEyessubnumber*\IEEEeqnarraynumspace \label{Eq: Y_i m N2N}\\
    \Y_{2,i} & = & c_{2,i}{\left( \X_i \right)^2},
    & \quad \forall i\IEEEeqnarraynumspace \label{Eq: Y_i p N2N}\\
    \X_i &=& \Om \cdot A_i,
    & \quad \forall i.
\end{IEEEeqnarray}
\noindent Let $\Expected{\Om} {=} M$, and since for a random variable $\boldsymbol{Z}$, $\Expected{\boldsymbol{Z}^2} {=} \text{Var}\left(\boldsymbol{Z}\right) + \Expected{\boldsymbol{Z}}^2$:
\begin{IEEEeqnarray}{rCll}
    \Expected{\Y_{1,i}} &=& -\left( 2 c_{2,i} p_i + c_{1,i}\right)\left(M \cdot A_i\right), &\forall i
    \IEEEyesnumber\IEEEyessubnumber*\\
    \Expected{\Y_{2,i}} &=& c_{2,i} \left(
    \left(M \cdot A_i\right)^2
    {+}
    \Var{\Om \cdot A_i} \right), &\forall i. \IEEEeqnarraynumspace \label{Eq: Exp Y2}
\end{IEEEeqnarray}
\noindent The variance of a linear combination can be calculated by
\begin{IEEEeqnarray}{rclll}
    \IEEEeqnarraymulticol{5}{C}{\Var{\sum_n \left( b_n \boldsymbol{Z}_n\right)} = \sum_{n,m}^N \left[ b_nb_m \text{cov}{\left(\boldsymbol{Z}_n\boldsymbol{Z}_m\right)}\right],} \label{Eq: Var LinearComb}\\
    \noalign{\noindent where $b_n$ is a scalar and $\boldsymbol{Z}_n$ a random variable.
    Let $\Sigma {=} \text{Var}\left[{\Om}\right]$ be the non-trivial covariance matrix of the forecast errors, we can then rewrite \eqref{Eq: Exp Y2} using \eqref{Eq: Var LinearComb} as
    \vspace{2\jot}}
    \Expected{\Y_{2,i}} &=& c_{2,i} \left[
    \left(M \cdot A_i\right)^2 {+}\AZsq\right], &~ \forall i. \IEEEeqnarraynumspace \label{Eq: Y2 Covar}
\end{IEEEeqnarray}
\noindent Finally, by using expression \eqref{Eq: Exp AB Yp Ym N2N} the expected costs are 
\begin{IEEEeqnarray}{rll}
    \Expected{c_i{(\boldsymbol{p}_{i})}} = c_i\left(p_{i}\right) 
    &-  {\left(
    2 c_{2,i} p_i  
    {+} c_{1,i}
    \right)}{\left( M \cdot A_i \right)} \IEEEnonumber \\
    &+ c_{2,i} {\left[ \left(M {\cdot} A_i\right)^2 {+} \AZsq
    \right]}&, \forall i. \IEEEeqnarraynumspace
\end{IEEEeqnarray}
%


\end{document}